\def\no{\if01}
\def\iftwelvept{\no}

\def\ifusepdf{\no}
\def\ifpsfont{\no}

\iftwelvept
\documentclass[12pt]{amsart}
\else
\documentclass[12pt]{amsart}
\fi
\usepackage{amssymb}

\usepackage{amscd}
\usepackage{latexsym}
\usepackage{verbatim}
\usepackage{mathrsfs}
\usepackage[all]{xy}

\ifusepdf
\usepackage{hyperref}
\else\fi
\ifpsfont
\usepackage[T1]{fontenc}
\usepackage{times}
\else\fi

\iftwelvept
\else\fi

\theoremstyle{plain}
\newtheorem{Theorem}{Theorem}[section]
\newtheorem{Proposition}[Theorem]{Proposition}
\newtheorem{Lemma}[Theorem]{Lemma}
\newtheorem{Corollary}[Theorem]{Corollary}

\newtheorem{Claim}{Claim}[Theorem]

\theoremstyle{definition}

\newtheorem{Definition}[Theorem]{Definition}
\newtheorem{Remark}[Theorem]{Remark}

\renewcommand{\theTheorem}{\arabic{section}.\arabic{Theorem}}

\newcommand{\ZZ}{{\mathbb{Z}}}
\newcommand{\QQ}{{\mathbb{Q}}}
\newcommand{\RR}{{\mathbb{R}}}

\newcommand{\NN}{{\mathbb{N}}}
\newcommand{\GG}{{\mathbb{G}}}

\newcommand{\nn}{{\mathfrak{n}}}
\newcommand{\pp}{{\mathfrak{p}}}
\newcommand{\XXX}{{\mathscr{X}}}

\newcommand{\VVV}{{\mathscr{V}}}
\newcommand{\DDD}{{\mathscr{D}}}

\newcommand{\OOO}{{\mathscr{O}}}
\newcommand{\OO}{{\mathcal{O}}}
\newcommand{\gp}{\textup{gp}}
\newcommand{\sm}{\textup{sm}}

\newcommand{\LLL}{\mathcal{L}}
\newcommand{\MMM}{\mathcal{M}}
\newcommand{\NNN}{\mathcal{N}}
\newcommand{\PPP}{\mathcal{P}}
\newcommand{\SSS}{\mathcal{S}}
\newcommand{\UU}{\mathcal{U}}

\newcommand{\Hom}{\textup{Hom}}

\newcommand{\Image}{\textup{Image}}

\newcommand{\Spec}{\operatorname{Spec}}

\newcommand{\mult}{\textup{mult}}

\newcommand{\rank}{\textup{rk}}

\newcommand{\XX}{\mathcal{X}}

\newcommand{\Proof}{{\sl Proof.}\quad}
\newcommand{\QED}{{\unskip\nobreak\hfil\penalty50\quad\null\nobreak\hfil
{$\Box$}\parfillskip0pt\finalhyphendemerits0\par\medskip}}

\begin{document}

\title[Log geometry and toric algebraic stacks]{Logarithmic geometry, minimal free resolutions \\ and toric algebraic stacks}


\author{Isamu Iwanari}





\begin{abstract}
In this paper we will introduce a certain type of morphisms of
log schemes (in the sense of Fontaine, Illusie and Kato)
and study their moduli.
Then by applying this we define the notion of toric algebraic stacks,
which may be regarded as torus and toroidal
emebeddings in the framework
of algebraic stacks and prove some fundamental properties.
Furthermore, we study the stack-theoretic analogue of toroidal embeddings.
\vspace{1mm}

\begin{flushleft}

\end{flushleft}

\end{abstract}

\maketitle

\section*{Introduction}
\renewcommand{\theTheorem}{\Alph{Theorem}}

In this paper we will introduce
a certain type of morphisms of log schemes
(in the sense of Fontaine, Illusie,
and Kato) and investigate their moduli.
Then by applying this we define a notion of {\it toric algebraic stacks}
over arbitrary schemes,
which may be regarded as {\it torus embeddings within the framework
of algebraic stacks}, and study some basic properties.
Our notion of toric algebraic stacks gives a natural generalization of
smooth torus embedding (toric varieties) over arbitrary schemes
preserving the smoothness,
and it is closely related to simplicial
toric varieties.
Moreover, our approach is amenable
and it also yields a sort of ``{\it stacky toroidal embeddings}".

\vspace{1mm}

We first introduce a 
notion of the {\it admissible and minimal free resolutions} of a monoid.
This notion plays a central role in this paper.
This leads to define a certain type of morphisms
of fine log schemes called {\it admissible FR} morphisms.
(``FR'' stands for {\it free resolution}.)
We then study the moduli stack of admissible FR
morphisms into a toroidal embedding endowed with the canonical log structure.
One may think of these moduli as a sort of natural ``stack-theoretic
generalization" of the classical notion of toroidal embeddings.

\vspace{1mm}

As promised above, the concepts of admissible FR
morphisms and their moduli
stacks yield the notion
 of {\it toric algebraic stacks} over arbitrary schemes.
Actually in the presented work on toric algebraic stacks,
admissible free resolutions of monoids
and admissible FR morphisms
play the role which is analogous to that of 
the monoids and monoid rings arising from cones
in classical toric geometry.
That is to say,
the algebraic aspect of toric algebraic stacks
is the algebra of admissible free resolutions of monoids.
In a sense, our notion of toric algebraic stacks is 
a hybrid of the definition of toric varieties given in \cite{AMRT}
and the moduli stack of admissible FR morphisms.
Fix a base scheme $S$.
Given a simplicial fan $\Sigma$
with additional data called ``level" $\nn$,
we define a toric algebraic stack $\XXX_{(\Sigma,\Sigma_{\nn}^0)}$.
It turns out that this stack has fairly good properties.
It is a {\it smooth Artin stack} of finite type
over $S$ with finite diagonal,
whose coarse moduli space is the simplicial toric variety
$X_{\Sigma}$ over $S$ (see Theorem~\ref{maintoric}).
Moreover it has a torus-embedding, a torus action functor
and a natural coarse moduli map, which
are defined in canonical fashions.
The complement of the torus
is a divisor with normal crossings relative to $S$.
If $\Sigma$ is non-singular and $\nn$ is a canonical level,
then $\XXX_{(\Sigma,\Sigma_{\nn}^0)}$ is the smooth toric variety $X_{\Sigma}$
over $S$.
Thus we obtain the following diagram of (2)-categories,
\[
\xymatrix@R=1mm @C=12mm{
  & (\textup{Toric algebraic stacks over}\ S)\ar[dd]^c \\
(\textup{Smooth toric varieties over}\ S)\ar[ur]^a \ar[dr]^b & \\
 &  (\textup{Simplicial toric varieties over}\ S) \\
}
\]
where $a$ and $b$ are fully faithful functors and
$c$ is an essentially surjective functor (see Remark~\ref{diag}).
One remarkable point to notice
is that working in the framework of algebraic stacks (including Artin stacks)
allows one to have a generalization of smooth toric varieties over $S$
that preserves the important features of smooth toric varieties
such as the smoothness.
There is another point to note.
Unlike toric varieties,
some properties of toric algebraic stacks depend very much
on the choice of a base scheme.
For example, the question of whether or not $\XXX_{(\Sigma,\Sigma_{\nn}^0)}$ is Deligne-Mumford depends on the base scheme.
Thus it is natural to develop our theory over arbitrary schemes.

\vspace{1mm}

Over the complex number field (and algebraically closed fields
of characteristic zero),
one can construct simplicial toric varieties
as geometric quotients by means of homogeneous coordinate
rings (\cite{C}).
In \cite{BCS}, by generalizing Cox's construction,
toric Deligne-Mumford stacks was introduced,
whose theory comes from Cox's viewpoint of toric varieties.
On the other hand,
roughly speaking,
our construction stemmed from the usual
definition of toric varieties given in, for example
\cite{KKMS}, \cite{AMRT}, \cite{F}, \cite{Oda}, \cite[Chapter IV, 2]{CF},
in log-algebraic geometry,
and it also yields a sort of ``stacky toroidal embeddings".
We hope that toric algebraic stacks provide
an ideal testing ground for problems and conjectures
on stacks in many areas of mathematics,
such as arithmetic geometry, algebraic geometry, mathematical physics,
etc.

\vspace{1mm}

This paper is organized as follows.
In section 2, we define the notion of
adimissible and minimal free resolution of monoids
and admissible FR morphisms
and investigate their properties.
It is an ``algebra" part of this paper.
In section 3, we construct algebraic moduli stacks
of admissible FR morphisms of toroidal embeddings
with canonical log structures.
In section 4, we define the notion of toric algebraic stacks
and prove fundamental properties by applying section 2 and 3.

\vspace{1mm}

{\it Applications and further works}.
Let us mention applications and further works,
which are not discussed in this paper.
The presented paper plays a central role in the subsequent papers
(\cite{I2}, \cite{I3}).
In \cite{I2}, by using the results and machinery
presented in this paper,
we study the 2-category of toric algebraic stacks,
and show that 2-category of toric algebraic stacks are equivalent to
the category of stacky fans.
Furthermore we prove that
toric algebraic stacks defined in this paper have a quite nice
geometric characterization in characteristic zero.
In \cite{I3},
we calculate the integral Chow ring of a toric algebraic stack
and show that it is isomorphic to the Stanley-Reisner ring.
As a possible application, we hope that
our theory might be applied to smooth toroidal
compactifications of spaces (including arithmetic schemes)
that can not be smoothly compactified in classical toroidal geometry
(cf. \cite{AMRT}, \cite{F}, \cite{N}).

\vspace{1mm}


{\bf Notations And Conventions}

\smallskip
\begin{enumerate}
\renewcommand{\labelenumi}{(\theenumi)}

\item We will denote by $\NN$ the set of {\it natural numbers}, by which
we mean the set of integers $n \ge 0$, by $\ZZ$
the {\it ring of rational integers}, by $\QQ$ the {\it rational number field},
and by $\RR$ the {\it real numbers}.
We write $\rank (L)$ for the rank of a free abelian group $L$.

\item By an {\it algebraic stack} we mean an algebraic stack
in the sense of \cite[4.1]{LM}.
All schemes, algebraic spaces, and algebraic stacks
are assumed to be {\it quasi-separated}.
We call an algebraic stack $\XXX$
 which admits an \'etale surjective cover
$X\to \XXX$, where $X$ is a scheme a {\it Deligne-Mumford stack}.
For details on algebraic stacks, we refer to \cite{LM}.
Let us recall the definition
of coarse moduli spaces and the fundamental existence theorem due to
Keel and Mori (\cite{KM}). Let $\XXX$ be an algebraic stack.
A {\it coarse moduli map (or space)} for $\XXX$ is
a morphism $\pi:\XXX\to X$ from $\XXX$ to an algebraic space $X$
such that the following conditions hold.
\begin{enumerate}
\item If $K$ is an algebraically closed field,
then the map $\pi$ induces a bijection between
the set of isomorphism classes of objects in $\XXX(K)$ and $X(K)$.

\item The map $\pi$ is universal for maps from $\XXX$ to algebraic spaces.
\end{enumerate}
Let $\mathcal{X}$ be an algebraic stack of finite type
over a locally noetherian scheme $S$ with finite diagonal.
Then a result of Keel and Mori says that
there exists a coarse moduli space
$\pi:\mathcal{X}\to X$ with $X$ of finite type and separated over $S$
(See also \cite{Con} in which the Noetherian assumption is eliminated).
Moreover $\pi$ is proper, quasi-finite and surjective, and the
natural map $\OO_X\to \pi_*\OO_{\mathcal{X}}$ is an isomorphism.
If $S'\to S$ is a flat morphism, then $\mathcal{X}\times_SS' \to S'$
is also a coarse moduli map.
\end{enumerate}

{\it Acknowledgement}.
I would like to thank Masao Aoki, Quo-qui Ito and Fumiharu Kato
for their valuable comments.
I also want to thank
Tadao Oda for his interest and offering me many excellent articles on
toric geometry.
I would like to thank
Institut de Math\'ematiques de Jussieu
for the hospitality during the stay where
a certain part of this work was done.
I am supported by Grant-in-Aid for JSPS Fellowships.

\renewcommand{\theTheorem}{\arabic{section}.\arabic{Theorem}}
\renewcommand{\thesubsubsection}{\arabic{section}.\arabic{subsection}.\arabic{subsubsection}}


\section{Toroidal and Logarithmic Geometry}
We first review  
some definitions and basic facts
concerning toroidal geometry and
logarithmic geometry in the sense of Fontaine, Illusie and K. Kato,
and establish notation for them.
We refer to \cite[Chapter IV. 2]{CF} \cite{F} \cite{KKMS} for details on toric and
toroidal geometry,
and refer to \cite{log} \cite{log2} for details on logarithmic geometry.

\subsection{Toric varieties over a scheme}
\label{TT}
Let $N\cong \ZZ^d$ be a lattice and $M:=\Hom_{\ZZ} (N,\ZZ)$ its dual.
Let $\langle \bullet,\bullet \rangle : M\times N\to \ZZ$ be the natural pairing.
Let $S$ be a scheme.
Let $\sigma\subset N_{\RR}:=N\otimes_{\ZZ}\RR$ be a strictly convex rational polyhedral cone
and 
\[
\sigma^{\vee}:=\{ m\in M_{\RR}:=M\otimes_{\ZZ}\RR \mid \langle m,u \rangle \ge
0\ \textup{ for all}\ u \in \sigma \}
\]
its dual.
(In this paper, all cones are assumed to be a
strictly convex rational polyhedral, unless otherwise stated.)
The {\it affine toric variety }(or {\it affine torus embedding}) $X_{\sigma}$
associated to $\sigma$ over $S$ is defined by
\[
X_{\sigma}:=\Spec\ \OO_S[\sigma^{\vee}\cap M]
\]
where $\OO_S[\sigma^{\vee}\cap M]$ is the monoid
algebra of $\sigma^{\vee} \cap M$
over the scheme $S$.

Let $\sigma\subset N_{\RR}$ be a cone.
(We sometimes use the $\QQ$-vector space $N\otimes_{\ZZ}\QQ$
instead of $N_{\RR}$.)
Let $v_1,\ldots,v_m$ be a minimal set of generators $\sigma$.
Each $v_i$ spans a {\it ray}, i.e., a 1-dimensional face of $\sigma$.
The affine toric variety $X_{\sigma}$ is smooth over $S$ if and
only if the first lattice points of $\RR_{\ge 0}v_1,\ldots,\RR_{\ge 0}v_m$
form a part of basis of $N$ (cf. \cite[page 29]{F}).
In this case, we refer to $\sigma$ as a {\it nonsingular cone}.
The cone $\sigma$ is {\it simplicial}
if it is generated by $\dim (\sigma)$ lattice points, i.e.,
$v_1,\ldots,v_m$ are linearly independent.
Let $\sigma$ be an $r$-dimensional simplicial cone in $N_{\RR}$
and $v_1,\ldots,v_r$ the first lattice points of rays in $\sigma$.
The {\it multiplicity} of $\sigma$, denoted by $\mult (\sigma)$, is defined to be 
the index $[N_{\sigma}:\ZZ v_1+\cdots +\ZZ v_r]$.
Here $N_{\sigma}$ is the lattice generated by $\sigma \cap N$.
If the multiplicity of a simplicial cone $\sigma$
is invertible on $S$, we say that the cone $\sigma$
is {\it tamely simplicial}.
If $\sigma$ and $\tau$ are cones, we write $\sigma\prec \tau$ (or $\tau \succ \sigma$)
to mean that $\sigma$ is a {\it face} of $\tau$.

A {\it fan} (resp. {\it simplicial fan, tamely simplicial fan}) $\Sigma$ in $N_{\RR}$ is a set of cones (resp. {\it simplicial cones, tamely simplicial cones}) in $N_{\RR}$
such that:
\begin{enumerate}
\renewcommand{\labelenumi}{(\theenumi)}

\item Each face of a cone in $\Sigma$ is also a cone in $\Sigma$,

\item The intersection of two cones in $\Sigma$ is a face of each.

\end{enumerate}

If $\Sigma$ is a fan in $N_{\RR}$, we denote by $\Sigma(r)$ the set of $r$-dimensional
cones in $\Sigma$, and denote by $|\Sigma|$ the
support of $\Sigma$ in $N_{\RR}$, i.e., the union of cones in $\Sigma$.
(Note that the set $\Sigma$ is not necessarily finite.
Even in classical situations, infinite fans are important and
they arises in various contexts
such as constructions of degeneration of abelian varieties,
the construction of hyperbolic Inoue surfaces, etc.)
\vspace{2mm}

Let $\Sigma$ be a fan in $N_{\RR}$.
There is a natural patching
of affine toric varieties associated to cones in $\Sigma$,
and the patching
defines a scheme of finite type and separated over $S$. 
We denote by $X_{\Sigma}$ this scheme, and we refer it as the
{\it toric variety} (or {\it torus embedding}) associated to $\Sigma$.
A toric variety $X$ contains a split algebraic torus
$T=\mathbb{G}_{m,S}^n=\Spec \OO_S[M]$ as an open dense subset, and the action of $T$
on itself extends to an action on $X_{\Sigma}$.

For a cone $\tau\in\Sigma$
we define its associated torus-invariant closed subscheme $V(\tau)$ to be the union
\[
\bigcup_{\sigma\succ \tau} \Spec\ \OO_S[(\sigma^{\vee}\cap M)/(\sigma^{\vee}\cap \tau_0^{\vee}\cap M)]
\]
in $X_{\Sigma}$, where $\sigma$ runs through the cones which contains $\tau$ as a face
and
\[
\tau_0^{\vee}:=\{m\in \tau^{\vee}| \langle m,n\rangle>0\ \textup{for some}\ n\in \tau\}
\]
(Affine schemes on the right hand naturally patch together, and the symbol $/(\sigma^{\vee}\cap \tau_0^{\vee}\cap M)$ means the ideal generated by $\sigma^{\vee}\cap \tau_0^{\vee}\cap M$).  We have
\[
\Spec\ \OO_S[(\sigma^{\vee}\cap M)/(\sigma^{\vee}\cap \tau_0^{\vee}\cap M)]
=\Spec\ \OO_S[\sigma^{\vee}\cap \tau^{\perp}\cap M]\subset \Spec\ \OO_S[\sigma^{\vee}\cap M]
\]
and the split torus $\Spec \OO_S[\tau^{\perp}\cap M]$ is a dense open subset
of $\Spec \OO_S[\sigma^{\vee}\cap \tau^{\perp}\cap M]$,
where $\tau^{\perp}=\{ m\in M\otimes_{\ZZ}\RR |\ \langle m,n \rangle =0\ 
\textup{for any}\ n\in \tau\}$. (Notice that $\tau_0^{\vee}=\tau^{\vee}-\tau^{\perp}$.)
For a ray $\rho\in \Sigma(1)$
(resp. a cone $\sigma\in \Sigma$),
we shall call $V(\rho)$ (resp. $V(\sigma)$)
the {\it torus-invariant divisor} (resp. {\it torus-invariant cycle})
associated to $\rho$ (resp. $\sigma$).
The complement $X_{\Sigma}-T$ set-theoretically equals
the union $\cup_{\rho\in \Sigma(1)}V(\rho)$.
Set $Z_{\tau}=\Spec \OO_S[\tau^{\perp}\cap M]$.
Then we have a natural stratification
$X_{\Sigma}=\sqcup_{\tau\in\Sigma}Z_{\tau}$
and for each cone $\tau\in \Sigma$, the locally closed subscheme
$Z_{\tau}$ is a $T$-orbit.

\subsection{Toroidal embeddings}
\label{toroidal}
Let $X$ be a normal variety over a field $k$,
i.e.,
a geometrically
integral normal scheme of finite type and separated over $k$.
Let $U$ be
a smooth Zariski open set of $X$.
We say that a pair $(X,U)$ is a
{\it toroidal embedding}
(resp. {\it good toroidal embedding, tame toroidal embedding})
if for every closed point
$x$ in $X$ there exist an \'etale neighborhood $(W,x')$ of $x$,
 an affine toric variety 
(resp. an affine simplicial toric variety, an affine tamely simplicial toric variety)
$X_{\sigma}$ over $k$, and an \'etale morphism
\[
f:W\longrightarrow X_{\sigma}.
\]
such that $f^{-1}(T_{\sigma})=W\cap U$. Here 
$T_{\sigma}$ is the algebraic torus in $X_{\sigma}$.

\subsection{Logarithmic geometry}

First of all, we shall recall some generalities on monoids.
In this paper, all monoids will assume to be commutative with unit.
Given a monoid $P$, we denote by $P^{\gp}$ the Grothendieck group of $P$.
If $Q$ is a submonoid of $P$,
we write $P\to P/Q$ for the cokernel in the category of monoids.
Two elements $p, p'\in P$ have the same image in $P/Q$
if and only if there exist $q, q'\in Q$
such that $p+q=p'+q'$.
The cokernel $P/Q$ has a monoid structure in the natural manner.
A monoid $P$ is {\it finitely generated} if there exists a surjective map
$\NN^{r}\to P$ for some positive integer $r$.
A monoid $P$ is said to be {\it sharp} if whenever $p+q=0$ for $p,q\in P$,
then $p=q=0$.
We say that $P$ is {\it integral}
if the natural map $P\to P^{\gp}$ is injective.
A finitely generated and integral monoid
is said to be {\it fine}.
An integral monoid $P$ is {\it saturated}
if for every $p \in P^{\gp}$ such that
$np\in P$ for some $n >0$,
it follows that $p\in P$.
An integral monoid $P$ is said to be
{\it torsion free} if $P^{\gp}$ is a torsion free abelian group.
We remark that a fine, saturated and sharp monoid is torsion free.

Given a scheme $X$, a {\it prelog structure} on $X$ is
a sheaf of monoids $\MMM$ on the \'etale site of $X$ together with
a homomorphism of sheaves of monoids $h:\MMM\to \OO_X$,
where $\OO_X$ is viewed as a monoid under multiplication.
A prelog structure is a {\it log structure}
if the map $h^{-1}(\OO_X^*)\to \OO_X^*$ is an isomorphism.
We usually denote simply by $\MMM$ the log structure $(\MMM,h)$
and by $\overline{\MMM}$ the sheaf $\MMM/\OO_X^*$.
A morphism of prelog structures $(\MMM,h)\to(\MMM',h')$
is a map $\phi:\MMM\to\MMM'$ of sheaves of monoids
such that $h'\circ \phi =h$.

For a prelog structure $(\MMM,h)$ on $X$,
we define its {\it associated log structure}
$(\MMM^{a},h^a)$ to be the push-out of
\[
\begin{CD}
 h^{-1}(\OO_X^*)  @>>> \MMM  \\
@VVV                      \\
\OO_X^* 
\end{CD}
\]
in the category of sheaves of monoids on the \'etale site $X_{et}$.
This gives the left adjoint functor of
the natural inclusion functor 
\[
(\textup{log structures on}\ X) \to
(\textup{prelog structures on}\ X).
\]

We say that a log structure $\MMM$ is {\it fine} if
\'etale locally on $X$ there exists a fine monoid
and a map $P\to \MMM$ from the constant sheaf associated to $P$ such
that $P^a\to \MMM$ is an isomorphism.
A fine log scheme $(X,\MMM)$ is {\it saturated}
if each stalk of $\MMM$ is a saturated monoid.
We remark that if $\MMM$ is fine and saturated, then
each stalk of $\overline{\MMM}$ is fine and saturated.

A morphism of log schemes
$(X,\MMM)\to (Y,\NNN)$ is
a pair $(f,h)$ of a morphism of underlying schemes $f:X\to Y$
and a morphism of log structures $h:f^*\NNN\to \MMM$,
where $f^*\NNN$ is the log structure associated to the composite $f^{-1}\NNN\to f^{-1}\OO_Y\to \OO_X$.
A morphism $(f,h):(X,\MMM)\to (Y,\NNN)$
is said to be {\it strict} if $h$ is an isomorphism.

Let $P$ be a fine monoid.
Let $S$ be a scheme.
Set $X_P:=\Spec \OO_{S}[P]$.
The {\it canonical log structure} $\MMM_P$ on $X_P$
is the fine log structure induced by the inclusion map $P\to \OO_S[P]$.
Let $\Sigma$ be a fan in $N_{\RR}$ ($N\cong\ZZ^d$) and $X_{\Sigma}$
the associated toric variety over $S$.
Then we have an induced log structure $\MMM_{\Sigma}$
on $X_{\Sigma}$ by
gluing the log structures arising from the homomorphism
$\sigma^{\vee}\cap M \to \OO_S[\sigma^{\vee}\cap M]$
for each cone $\sigma\in \Sigma$. Here $M=\Hom_{\ZZ} (N,\ZZ)$.
We shall refer this log structure as the
{\it canonical log structure} on $X_{\Sigma}$.
If $S$ is a locally noetherian regular scheme, we have that $\MMM_{\Sigma}=\OO_{X_{\Sigma}}\cap i_*\OO_{\Spec \OO_{S}[M]}^*$ where $i:\Spec \OO_S[M]\to X_{\Sigma}$
is the torus embedding
(cf. \cite[11.6]{log2}).

Let $(X,U)$ be a toroidal embedding over a field $k$
and $i:U\to X$ be the natural immersion.
Define a log structure $\alpha_X:\MMM_X:=\OO_X\cap i_*\OO_U^* \to \OO_X$
on $X$.
This log structure is fine and saturated and said to be the
{\it canonical log structure} on $(X,U)$.


\section{Free Resolution of Monoids}

\subsection{Minimal and admissible free resolution of a monoid}

\begin{Definition}
Let $P$ be a monoid.
The monoid $P$ is said to be {\it toric} if $P$ is a fine, saturated and torsion free monoid. 

\end{Definition}

\begin{Remark}
\label{monoid1}
If a monoid $P$ is toric,
there exists a strictly convex rational polyhedral 
cone $\sigma \in \Hom_{\ZZ} (P^{\gp},\ZZ)\otimes_{\ZZ}\QQ$
such that $\sigma^{\vee} \cap P^{\gp} \cong P$.
Here the dual cone $\sigma^{\vee}$ lies on $P^{\gp}\otimes_{\ZZ}{\QQ}$.
Indeed, we see this as follows.
There exists a sequence of canonical injective homomorphisms
$P\to P^{\textup{gp}}\to P^{\textup{gp}}\otimes_{\ZZ}\QQ$.
Define a cone
\[
C(P):=\{ \Sigma_{i=0}^{n}a_i\cdot p_i |\ a_i\in \QQ_{\ge0},\ 
p_i\in P\}\subset P^{\textup{gp}}\otimes_{\ZZ}\QQ.
\]
Note that it is a {\it full-dimensional} rational polyhedral cone
(but not necessarily strictly convex),
and $P=C(P)\cap P^{\textup{gp}}$ since $P$ is saturated.
Thus the dual cone $C(P)^{\vee}\subset \Hom_{\ZZ} (P^{\gp},\ZZ)\otimes_{\ZZ}\QQ$
is a strictly convex rational polyhedral cone (cf. \cite[(13) on page 14]{F}).
Hence our assertion follows.
\end{Remark}

$\bullet$ Let $P$ be a monoid and $S$ a submonoid of $P$.
We say that the submonoid $S$ is {\it close to} the monoid $P$
if for every element
$e$ in $P$, there exists a positive integer $n$ such that
$n\cdot e$ lies in $S$.

$\bullet$ Let $P$ be a toric sharp monoid, and
let $r$ be the rank of $P^{\gp}$. 
A toric sharp monoid $P$ is said to be
{\it simplicially toric} if there
exists a submonoid $Q$ of $P$ generated by $r$ elements
such that $Q$ is close to $P$.

\begin{Lemma}
\label{simplicialtoric}
\begin{enumerate}
\renewcommand{\labelenumi}{(\theenumi)}
\item A toric sharp monoid $P$ is simplicially toric
if and only if
we can choose a (strictly convex rational polyhedral) 
simplicial full-dimensional cone 
\[
\sigma\subset \Hom_{\ZZ}(P^{\gp},\ZZ)\otimes_{\ZZ}\QQ
\]
such that
$\sigma^{\vee}\cap P^{\gp}\cong P$, where $\sigma^{\vee}$ denotes
the dual cone in $P^{\gp}\otimes_{\ZZ}\QQ$.

\item If $P$ is a simplicially toric sharp monoid, then
\[
C(P):=\{ \Sigma_{i=0}^{n}a_i\cdot p_i |\ a_i\in \QQ_{\ge0},\ 
p_i\in P\}\subset P^{\textup{gp}}\otimes_{\ZZ}\QQ
\]
is a  (strictly convex rational polyhedral) simplicial full-dimensional cone.
\end{enumerate}

\end{Lemma}

\Proof
We first prove (1).
The ``if" direction is clear.
Indeed, if there exists such a simplicial full-dimensional cone $\sigma$,
then the dual cone $\sigma^{\vee}$ is also a simplicial full-dimensional cone.
Let $Q$ be the submonoid of $P$, which is generated by
the first lattice points of rays on $\sigma^{\vee}$.
Then $Q$ is close to $P$.

Next we shall show the ``only if" part.
Assume there exists a submonoid $Q\subset P$
such that $Q$ is close to $P$ and generated by $\rank (P^{\textup{gp}})$
elements.
By Remark~\ref{monoid1}
there exists a cone 
$C(P)=\{ \Sigma_{i=0}^{n}a_i\cdot p_i |\ a_i\in \QQ_{\ge0},\ 
p_i\in P\}$
such that
$P=C(P)\cap P^{\textup{gp}}$.
Note that since $P$ is sharp,
$C(P)$ is strictly convex and full-dimensional.
Thus $\sigma:=C(P)^{\vee}\subset \Hom_{\ZZ}(P^{\gp},\ZZ)\otimes_{\ZZ}\QQ$ is a
full-dimensional cone.
It suffices to show that $C(P)$ is simplicial, i.e.,
the cardinality of the set of rays of $C(P)$ is equal to
the rank of $P^{\gp}$.
For any ray $\rho$ of $\sigma^{\vee}=C(P)$, $Q\cap \rho$ is non-empty because
$Q$ is close to $P$.
Thus $Q$ can not be generated
by any set of elements of $Q$ whose cardinality is less than the cardinality
of rays in $\sigma^{\vee}$.
Thus we have $\rank (P^{\gp})\ge \#\sigma^{\vee}(1)$
(here we write $\#\sigma^{\vee}(1)$ for the cardinality of the set of rays of $\sigma^{\vee}$).
Hence $\sigma^{\vee}=C(P)$ is simplicial and thus
$\sigma$ is also simplicial. It follows (1).
The assertion of (2) is clear.
\QED

\begin{Lemma}
\label{rankremark}
Let $P$ be a toric sharp monoid. Let $F$ be a monoid such that
$F\cong \NN^r$ for some $r\in \NN$.
Let $\iota:P\to F$ be an injective homomorphism
such that $\iota(P)$ is close to $F$.
Then the rank of $P^{\gp}$ is equal to the rank of $F$, i.e.,
$\rank (F^{\gp})=r$.
\end{Lemma}

\Proof
Note first that
$P^{\gp}\to F^{\gp}$ is injective.
Indeed, the natural homomorphisms $P\to F$
and $F\to F^{\gp}$ are injective.
Thus if $p_1, p_2\in P$ have
the same image in $F^{\gp}$, then $p_1=p_2$.
Hence $P^{\gp}\to F^{\gp}$ is injective.
Since $i(P)$ is close to $F$,
the cokernel of $P^{\gp}\to F^{\gp}$
is finite.
Hence our claim follows.
\QED

\begin{Proposition}
\label{MFR}
Let $P$ be a simplicially toric sharp monoid.
Then there exists
an injective homomorphism of monoids
\[
i:P \longrightarrow F
\]
which has the
following properties.
\begin{enumerate}
\renewcommand{\labelenumi}{(\theenumi)}

\item The monoid $F$ is isomorphic to $\NN^d$ for some $d\in \NN$,
and the submonoid $i(P)$ is close to $F$.

\item If $j:P \to G$ is an injective homomorphism,
and $G$ is isomorphic to $\NN^d$ for some $d\in \NN$, and $j(P)$ is  close to $G$,
then there exists a unique homomorphism $\phi:F\to G$
such that the diagram
\[
 \xymatrix{
 P \ar[d]_{j} \ar[r]^{i} & F \ar[dl]^{\phi}   \\
 G  \\
  }
\]
commutes.
\end{enumerate}

Furthermore if $C(P):=\{ \Sigma_{i=0}^{n}a_i\cdot p_i |\ a_i\in \QQ_{\ge0},\ 
p_i\in P\}\subset P^{\textup{gp}}\otimes_{\ZZ}\QQ$
$($it is a simplicial cone $(cf$. Lemma~\ref{simplicialtoric} $(2))$,
then there exists a canonical injective map
\[
F\to C(P)
\]
that has the following properties:
\begin{enumerate}
\renewcommand{\labelenumi}{(\alph{enumi})}

\item The natural diagram
\[
 \xymatrix{
 P \ar[d] \ar[r]^{i} & F \ar[dl]   \\
 C(P)  \\
  }
\]
commutes,
\item Each irreducible element of $F$ lies on a unique ray of $C(P)$
via $F\to C(P)$.
\end{enumerate}

\end{Proposition}
\Proof
Let $d$ be the rank of the torsion-free abelian group $P^{\gp}$.
By Lemma~\ref{simplicialtoric},
$C(P)$ is a full-dimensional simplicial cone in $P^{\gp}\otimes_{\ZZ}\QQ$.
Let $\{ \rho_1,\ldots,\rho_d\}$ be the set of rays in $C(P)$.
Let us denote by $v_i$ the first lattice point on $\rho_i$ in $C(P)$.
Then for any element $c\in C(P)$
we have a unique representation of $c$ such that
$c=\Sigma_{1\le i\le d}a_i\cdot v_{i}$ where $a_i\in \QQ_{\ge 0}$
for $1\le i\le d$.
Consider the map $q_k:P^{\gp}\otimes_{\ZZ}\QQ\to \QQ;\ p=\Sigma_{1\le i\le d}a_i\cdot v_{i}\mapsto a_{k}$ ($a_{i}\in \QQ$ for all $i$).
Set $P_i:=q_i(P^{\gp})\subset \QQ$.
It is a free abelian group generated by one element.
Let $p_i\in P_i$ be the element such that $p_i>0$ and the absolute value
of $p_i$ is the smallest in $P_i$.
Let $F$ be the monoid generated by $p_1\cdot v_1,\ldots,p_d\cdot v_d$.
Clearly,
we have $F\cong \NN^d$ and $P\subset F\subset P^{\gp}\otimes_{\ZZ}\QQ$.
Note that there exists a positive integer $b_i$
such that $p_i=1/b_i$ for each $1\le i\le d$.
Therefore $b_i\cdot p_i\cdot v_i=v_i$
for all $i$, thus it follows that $P$ is close to $F$.

It remains to show that $P\subset F$ satisfies the property (2).
Let $j:P\to G$ be an injective homomorphism of monoids such that
$j(P)$ is close to $G$.
Notice that by Lemma~\ref{rankremark}, we have $G\cong \NN^d$.
The monoid $P$ has the natural injection
$\iota:P\to P^{\gp}\otimes_{\ZZ}{\QQ}$.
On the other hand, for any element $e$ in $G$, there exists a positive integer
$n$ such that $n\cdot e$ is in $j(P)$. Therefore we have a unique injective homomorphism
 $\lambda:G\to P^{\gp}\otimes_{\ZZ}\QQ$ which extends $P\to P^{\gp}\otimes_{\ZZ}{\QQ}$ to
$G$. Indeed, if $g\in G$ and $n\in \ZZ_{\ge 1}$ such that $n\cdot g\in j(P)=P$, then
we define $\lambda(g)$ to be $\iota(n\cdot g)/n$
(it is easy to see that $\lambda(g)$ does not depend on
the choice of $n$).
The map $\lambda$ defines a homomorphism of monoids.
Indeed, if $g_1, g_2\in G$, then there exists a positive integer
$n$ such that both $n\cdot g_1$ and $n\cdot g_2$ lie in $P$,
and it follows that $\lambda(g_1+g_2)=\iota(n(g_1+g_2))/n=\iota(n\cdot g_1)/n+\iota(n\cdot g_2)/n=\lambda(g_1)+\lambda(g_2)$.
In addition, $\lambda$ sends the unit element of $G$ to the unit element
of $P^{\gp}\otimes_{\ZZ}{\QQ}$.
Since $P^{\gp}\otimes_{\ZZ}\QQ\cong\QQ^{d}$, thus
$\lambda:G\to P^{\gp}\otimes_{\ZZ}\QQ$ is a unique extension
of the homomorphism $\iota:P\to P^{\gp}\otimes_{\ZZ}{\QQ}$.
The injectiveness of $\lambda$ follows from its definition.
We claim that there exists a sequence of inclusions
\[
P \subset F \subset G \subset P^{\gp}\otimes_{\ZZ}\QQ.
\]
Since $P$ is close to $G$ and $C(P)$ is a full-dimensional
simplicial cone, thus each irreducible element of $G$
lies in a unique ray of $C(P)$.
(For a ray $\rho$ of $C(P)$, the first point of
$G\cap \rho$ is an irreducible element of $G$.)
On the other hand, we have
$\ZZ_{\ge 0}\cdot p_i\subset q_{i}(G^{\gp})\cap \QQ_{\ge 0}$.
This implies $F\subset G$. Thus we have (2).

By the above construction, clearly 
there exists the natural homomorphism $F\to C(P)$.
The property (a) is clear.
The property (b) follows from the above argument.
Hence we complete the proof of our Proposition.
\QED

\begin{Definition}
Let $P$ be a simplicially toric sharp monoid.
If an injective homomorphism of monoids
\[
i:P\longrightarrow F
\]
that
satisfies the properties $(1)$, $(2)$ (resp. the property $(1)$) in Proposition~\ref{MFR},
we say that $i:P\longrightarrow F$
is a {\it minimal free resolution} (resp. {\it admissible free resolution})
of $P$.
\end{Definition}

\begin{Remark}
\label{MFRR}

\begin{enumerate}
\renewcommand{\labelenumi}{(\theenumi)}

\item By the observation in the proof of Proposition~\ref{MFR},
if $j:P\to G$ is an admissible free resolution of a simplicially toric sharp
monoid $P$, then there is a natural commutative diagram
\[
 \xymatrix{
 P \ar[d] \ar[r]^{i} & F \ar[r]^{\phi} \ar[dl] & G \ar[dll]  \\
 C(P)  \\
  }
\]
such that $\phi\circ i=j$, where $i:P\to F$ is the minimal free resolution
of $P$. Furthermore, all three maps into $C(P)$ are injective and
each irreducible element of $G$ lies on a unique ray of $C(P)$.

\item By Lemma~\ref{rankremark}, the rank of $F$ is equal to the rank of $P^{\gp}$.

\item We define the {\it multiplicity} of $P$, denoted by $\mult (P)$,
to be the order of the cokernel of $i^{\gp}:P^{\gp}\to F^{\gp}$.
If $P$ is isomorphic to $\sigma^{\vee}\cap M$ where $\sigma$ is a simplicial cone, it is easy to see that $\mult (P)=\mult (\sigma)$.

\end{enumerate}
\end{Remark}

\begin{Proposition}
\label{MFR3}
Let $P$ be a simplicially toric sharp monoid and 
$i:P\to F\cong \NN^d$ its minimal free resolution.
Consider the following diagram
\[
 \xymatrix{
 P \ar[d]_{q} \ar[r]^{i} & F\cong\NN^d \ar[d]^{\pi}   \\
 Q \ar[r]^{j} & \NN^r, \\
  }
\]
where $Q:=\Image(\pi \circ i)$ and $\pi:\NN^d\to\NN^r$ is defined by
$(a_1,\ldots,a_d) \mapsto (a_{\alpha(1)},\ldots,a_{\alpha(r)})$.
Here $\alpha(1),\ldots,\alpha(r)$ are positive integers such that $1\le \alpha(1) < \cdots <\alpha(r)\le d$.
Then $Q$ is a simplicially
toric sharp monoid and $j$ is the minimal free resolution of $Q$.
\end{Proposition}

\Proof
After reordering we assume that $\alpha(k)=k$ for $1\le k\le r$.
First, we will show that $Q$ is a simplicially
toric sharp monoid.
Since $Q$ is close to $\NN^r$ via $j$,
$Q$ is sharp.
If $e_{i}$ denotes the $i$-th standard irreducible element,
then for each $i$, there exists positive integers $n_1,\ldots,n_r$
such that $n_1\cdots e_1,\ldots,n_r\cdot e_r\in Q$
and $n_1\cdots e_1,\ldots,n_r\cdot e_r$ generates
a submonoid which is close to $Q$.
Thus it suffices only to prove that $Q$ is a toric monoid.
Clearly, $Q$ is a fine monoid.
To see the saturatedness
we first regard $P^{\gp}$ and $Q^{\gp}$ as subgroups
of
$\QQ^{d}=(\NN^{d})^{\gp}\otimes_{\ZZ}\QQ$ and $\QQ^r=(\NN^{r})^{\gp}\otimes_{\ZZ}\QQ$ respectively.
It suffices to show $Q^{\gp}\cap \QQ^r_{\ge 0}=Q$.
Since $P$ is saturated, thus $P=P^{\gp}\cap \QQ^{d}_{\ge 0}$.
Note that $q^{\gp}:P^{\gp}\to Q^{\gp}$ is surjective.
It follows that $P^{\gp}\cap \QQ^{d}_{\ge 0}\to Q^{\gp}\cap \QQ^r_{\ge 0}$
is surjective.
Indeed, let $\xi\in  Q^{\gp}\cap \QQ^r_{\ge 0}$ and
$\xi'\in P^{\gp}$ such that $\xi=q^{\gp}(\xi')$.
Put $\xi'=(b_1,\ldots,b_d)\in P^{\gp}\subset (\NN^d)^{\gp}=\ZZ^{d}$.
Note that $b_i\ge 0$ for $1\le i\le r$.
Since $P^{\gp}$ is a subgroup of $\ZZ^d$ of a finite index,
there exists an element $\xi''=(0,\ldots,0, c_{r+1},\ldots,c_d)\in P^{\gp}$
such that
$\xi'+\xi''=(b_1,\ldots,b_{r},b_{r+1}+c_{r+1},\ldots,b_{d}+c_{d})\in \ZZ^d_{\ge 0}$.
Then $\xi=q^{\gp}(\xi')=q^{\gp}(\xi'+\xi'')$.
Thus $P^{\gp}\cap \QQ^{d}_{\ge 0}\to Q^{\gp}\cap \QQ^r_{\ge 0}$
is surjective.
Hence $Q$ is saturated.

It remains to prove that $Q\subset \NN^r$ is the minimal free resolution.
It order to prove this, recall that the construction of minimal free resolution of $P$.
With the same notation as in the first paragraph of the proof of Proposition~\ref{MFR}, the monoid $F$ is defined to be a free submonoid
$\NN \cdot p_1\cdot v_1\oplus \cdots \oplus\NN\cdot p_d\cdot v_d$ of $C(P)$
where $p_k\cdot v_k$ is the first point of $C(P)\cap \tilde{q}_k(P^{\gp})$ for $1\le k \le d$.
Here the map $\tilde{q}_k:P^{\gp}\otimes_{\ZZ}\QQ\to \QQ\cdot v_k$
is defined by $\Sigma_{1\le i\le d}a_i\cdot v_{i}\mapsto a_{k}\cdot v_k$ ($a_{i}\in \QQ$ for all $i$).
We shall refer this construction as the {\it canonical construction}.
After reordering, we have the following diagram
\[
 \xymatrix{
 \NN \cdot p_1\cdot v_1\oplus \cdots \oplus\NN\cdot p_d\cdot v_d\ \ar[d]_{\pi} \ar[r] & C(P) \ar[r] \ar[d]  & P^{\gp}\otimes_{\ZZ}\QQ \ar[d]\\
 \NN \cdot p_1\cdot v_1\oplus \cdots \oplus\NN\cdot p_r\cdot v_r \ar[r] & C(Q) \ar[r] &    Q^{\gp}\otimes_{\ZZ}\QQ \\
  }
\]
where $p_k\cdot v_k$ is regarded as a point on a ray of $C(Q)$
for $1\le i\le r$.
Then $p_k\cdot v_k$ is the first point of $C(Q)\cap \tilde{q}'_k(Q^{\gp})$ for $1\le k \le r$, where $\tilde{q}'_k:Q^{\gp}\otimes_{\ZZ}\QQ\to \QQ\cdot v_k;\ \Sigma_{1\le i\le r}a_i\cdot v_{i}\mapsto a_{k}\cdot v_k$ (note that
for any $c\in Q^{\gp}\otimes_{\ZZ}\QQ$, there is a unique representation of $c$ such that
$c=\Sigma_{1\le i\le r}a_i\cdot v_{i}$ where $a_i\in \QQ_{\ge 0}$
for $1\le i\le r$.). Then $Q\to \NN^r\cong \NN \cdot p_1\cdot v_1\oplus \cdots \oplus\NN\cdot p_r\cdot v_r$ is
the canonical construction for $Q$, and thus
it is the minimal free resolution. Hence we obtain our Proposition.
\QED

\begin{Proposition}
\label{AFR}

\begin{enumerate}
\renewcommand{\labelenumi}{(\theenumi)}

\item Let $\iota:P\to F$ be an admissible free resolution.
Then $\iota$ has the form
\[
\iota=n\circ i:P\stackrel{i}{\to} F\cong \NN^{d}\stackrel{n}{\to} \NN^{d}\cong F
\]
where $i$ is the minimal free resolution and $n:\NN^d\to\NN^d$
is defined by $e_{i}\mapsto n_i\cdot e_{i}$.
Here $e_i$ is the $i$-th standard irreducible element of $\NN^d$
and $n_i\in \ZZ_{\ge 1}$ for $1\le i\le d$.

\item Let $\sigma$ be a full-dimensional
simplicial cone in $N_{\RR}$ ($N=\ZZ^d$, $M=\Hom_{\ZZ}(N,\ZZ)$)
and $\sigma^{\vee}\cap M\hookrightarrow F$ the minimal free resolution
(note that $\sigma^{\vee}\cap M$ is a simplicially toric sharp monoid).
Then
there is a natural inclusion
$\sigma^{\vee}\cap M\subset F\subset \sigma^{\vee}$.
Each irreducible element of $F$ lies on a unique ray
of $\sigma^{\vee}$. This gives a bijective map
between the set of irreducible elements of $F$
and the set of rays of $\sigma^{\vee}$.

\end{enumerate}

\end{Proposition}

\Proof
We first show (1).
By Remark~\ref{monoid1} (1), there exist natural inclusions
\[
P\subset F \subset F\subset C(P)
\]
where the first inclusion $P\subset F$ is the minimal free resolution
and the composite $P\subset F\subset F$ is equal to $\iota:P\to F$.
Moreover each irreducible element of the left $F$
(resp. the right $F$) lies on a unique ray of $C(P)$.
Let $\{ s_1,\ldots ,s_d\}$ (resp. $\{ t_1,\ldots, t_d\}$)
denote images of irreducible elements of the left $F$ (resp. the right $F$)
in $C(P)$.
Since the rank of the free monoid $F$ is equal to the
cardinality of the set of rays of $C(P)$,
thus there is a positive integer $n_i$ such that
$n_i\cdot t_i\in \{s_1,\ldots ,s_d\}$ for $1\le i\le d$.
After reordering, we have $n_i\cdot t_i=s_i$ for each $1\le i\le d$.
Therefore our assertion follows.

To see (2), consider
\[
\sigma^{\vee}\cap M\subset C(\sigma^{\vee}\cap M)\subset (\sigma^{\vee}\cap M)^{\gp}\otimes_{\ZZ}\QQ =M\otimes_{\ZZ}\QQ\subset M\otimes_{\ZZ}\RR
\]
where $C(\sigma^{\vee}\cap M)=\{ \Sigma_{1\le i\le m} a_i\cdot s_i |\ a_i\in\QQ_{\ge 0},\ s_i\in \sigma^{\vee}\cap M\}$ is a simplicial
full-dimensional cone by Lemma~\ref{simplicialtoric}. The cone $\sigma^{\vee}$
is the completion of $C(\sigma^{\vee}\cap M)$ with respect to
the usual topology on $M\otimes_{\ZZ}\RR$.
Then the second assertion follows from Proposition~\ref{MFR} (b)
and the fact that
the rank of the free monoid $F$ is equal to the
cardinality of the set of rays of $C(\sigma^{\vee}\cap M)$.
\QED

Let $P$ be a toric monoid.
Let $I\subset P$ be an {\it ideal}, i.e.,
a subset such that $P+I\subset I$.
We say that $I$ is a {\it prime ideal} if $P-I$ is a submonoid of $P$.
Note that the empty set is a prime ideal.
Set $V=P^{\gp}\otimes_{\ZZ}\QQ$. 
For a subset $S\subset V$,
let $C(S)$ be the (not necessarily strictly convex) cone define by
$C(S):=\{ \Sigma_{1\le i\le n} a_i\cdot s_i |\ a_i\in\QQ_{\ge 0},\ s_i\in S\}$.
To a prime ideal $\pp\subset P$ we associate
$C(P-\pp)$.
By an elementary observation, we see that
$C(P-\pp)$ is a face of $C(P)$ and
it gives rise to a bijective correspondence
between the set of prime ideals of $P$ and the set of faces
of $C(P)$  (cf. \cite[Proposition 1.10]{Th}).

Let $P$ be a simplicially toric sharp monoid.
Then the cone $C(P)\subset P^{\gp}\otimes_{\ZZ}\QQ$
is a strictly convex rational polyhedral simplicial full-dimensional cone (cf. Lemma~\ref{simplicialtoric}).
A prime ideal $\pp\in P$ is called a {\it height-one prime ideal}
if $C(P-\pp)$ is a $(\dim P^{\gp}\otimes_{\ZZ}\QQ-1)$-dimensional face
of $C(P)$, equivalently $\pp$ is a minimal nonempty prime.
In this case,
for each height-one prime ideal $\pp$ of $P$
there exists a unique ray of $C(P)$, which does not
lie in $C(P-\pp)$.
Let $P\to F$ be the minimal free resolution.
Notice that the rank of $F$ is equal to the cardinality
of the set of rays of $C(P)$.
Therefore taking account of Proposition~\ref{MFR} (b),
there exists a natural bijective correspondence
between the set of rays of $C(P)$
and the set of irreducible elements of $F$.
Therefore
there exists the natural correspondences
\[
 \xymatrix@R=6mm @C=12mm{
\{ \textup{The set of height-one prime ideals of }P \}  \ar[d]^{\cong} \\
 \{ \textup{The set of rays of }C(P)\} \ar[d]^{\cong}\\
  \{ \textup{The set of irreducible elements of }F\}.\\
  }
\]

\begin{Definition}
\label{admtype}
Let $P$ be a simplicially toric sharp monoid.
Let $I$ be the set of height-one prime ideals of $P$.
Let $j:P\to F$
be an admissible free resolution of $P$.
Let us denote by $e_i$
the irreducible element of $F$
corresponding to $i\in I$.
We say that
$j:P\to F$ is an {\it admissible
free resolution of type $\{n_{i}\in\ZZ_{\ge 1}\}_{i\in I}$}
if $j$ is isomorphic to the composite $P\stackrel{i}{\to}F\stackrel{w}{\to}F$
where $i$ is the minimal free resolution and $w:F\to F$
is defined by $e_i\mapsto n_{i}\cdot e_i$.

Note that admissible free resolutions of a simplicially toric sharp monoid
are classified by their {\it type}.

\end{Definition}

We use the following technical Lemma in the subsequent section.

\begin{Lemma}
\label{monogp}
Let $P$ be a toric monoid and $Q$ a saturated subomonoid that is close to $P$.
Then the monoid $P/Q$ (cf. section 1.3)
is an abelian group, and
the natural homomorphism $P/Q\to P^{\gp}/Q^{\gp}$ is an isomorphism.
\end{Lemma}

\Proof
Clearly, $P/Q$ is finite and thus it is an abelian group.
We will prove that $P/Q\to P^{\gp}/Q^{\gp}$ is injective.
It suffices to show that $P\cap Q^{\gp}=Q$
in $P^{\gp}$.
Since $P\cap Q^{\gp}\supset Q$,
we will show $P\cap Q^{\gp}\subset Q$.
For any $p\in P\cap Q^{\gp}$, there exists a positive integer $n$
such that $n\cdot p\in Q$ because $Q$ is close to $P$.
Since $Q$ is saturated, we have $p\in Q$. Hence $P\cap Q^{\gp}=Q$.
Next we will prove that $P/Q\to P^{\gp}/Q^{\gp}$
is surjective. Let $p\in P^{\gp}$.
Take $p_1,p_2\in P$ such that $p=p_1-p_2$ in $P^{\gp}$.
It is enough to show that there exists $p'\in P$
such that $p'+p_2\in Q$.
Since $Q$ is close to $P$, thus our assertion is clear.
Hence $P/Q\to P^{\gp}/Q^{\gp}$ is sujective.
\QED


\subsection{MFR morphisms and admissible FR morphisms}
The notions defined below play a pivotal role in our theory.

\begin{Definition}
Let $(F,\Phi):(X,\mathcal{M})\to (Y,\mathcal{N})$ be a morphism 
of fine log-schemes. We say that $(F,\Phi)$ is an MFR (=Minimal Free
Resolution) morphism
if for any point $x$ in $X$, the monoid
$F^{-1}\overline{\mathcal{N}}_{\bar{x}}$ is simplicially toric and the homomorphism
of monoids $\overline{\Phi}_{\bar{x}}: F^{-1}\overline{\mathcal{N}}_{\bar{x}}\to \overline{\mathcal{M}}_{\bar{x}}$ is the minimal free resolution
 of $F^{-1}\overline{\mathcal{N}}_{\bar{x}}$.
\end{Definition}
\begin{Proposition}
\label{ex1}
Let $P$ be a simplicially toric sharp monoid.
Let $i:P\to F$ be its minimal free resolution.
Let $R$ be a ring.
Then the map $i:P\to F$ defines an MFR morphism of fine log schemes
$(f,h):(\Spec R[F],\MMM_F)\to (\Spec R[P],\MMM_P)$,
where $\MMM_F$ and $\MMM_P$ are log structures
induced by charts $F\to R[F]$ and $P\to R[P]$
respectively.
\end{Proposition}

\Proof
Since $\MMM_F$ and $\MMM_P$ are Zariski log structures
arising from $F\to R[F]$ and $P\to R[P]$
respectively, to prove our claim it suffices to consider
only Zariski stalks of log structures
i.e., to show that for any point of $x\in \Spec R[F]$
the homomorphism $h:f^{-1}\overline{\MMM}_{P,f(x)}\to \overline{\MMM}_{F,x}$
is the minimal free resolution.
Suppose that $F=\NN^r\oplus \NN^{d-r}$
and
$x \in\Spec R[(\NN^{d-r})^{\textup{gp}}]\subset \Spec R[\NN^{d-r}]=\Spec R[\NN^r\oplus \NN^{d-r}]/(\NN^{r}-\{0\})\subset \Spec R[\NN^r\oplus \NN^{d-r}]$.
Then $f(x)$ lies in
$\Spec R[P_0^{\textup{gp}}]\subset \Spec R[P_0]=\Spec R[P]/(P_1)\subset \Spec
R[P]$, where $P_0$ is the submonoid of elements whose images
of $u:P\to F=\NN^r\oplus\NN^{d-r}\stackrel{\textup{pr}_1}{\to}\NN^{r}$
are zero
and $P_1$ is the ideal generated by elements of $P$, whose images of
$u$
are non-zero.
Indeed, since $x \in\Spec R[\NN^r\oplus \NN^{d-r}]/(\NN^{r}-\{0\})$,
thus $f(x)\in \Spec R[P]/(P_1)$. For any $p\in P_0$, the
image of $i(p)$ in $(\NN^{d-r})^{\gp}$ is invertible, and thus $f(x)\in \Spec R[P_0^{\gp}]$.
Note that there exists the commutative diagram
\[
 \xymatrix{
 P \ar[d] \ar[r]^(0.3){i} & F=\NN^r\oplus \NN^{d-r} \ar[d]^{\operatorname{pr}_1}  \\
 f^{-1}\overline{\MMM}_{P,f(x)} \ar[r]^{\overline{h}} & \overline{\MMM}_{F,x}=\NN^r, \\
  }
\]
where the vertical surjective homomorphisms are
induced by the standard charts $P\to \MMM_{P}$ and $F\to \MMM_F$ respectively.
Applying Proposition~\ref{MFR3} to this diagram,
it suffices to prove that  
$f^{-1}\overline{\MMM}_{P,f(x)}\to\overline{\MMM}_{F,x}$ is injective.
Since there are a sequence of surjective maps
$P^{\gp}\to P^{\gp}/P_0^{\gp}\to f^{-1}\overline{\MMM}_{P,f(x)}^{\gp}$
and the inclusion $f^{-1}\overline{\MMM}_{P,f(x)}\subset f^{-1}\overline{\MMM}_{P,f(x)}^{\gp}$, thus
it is enough to prove that for any $p_1$ and $p_2$ in $P$
such that $u(p_1)=u(p_2)$ the element $p_1-p_2\in P^{\textup{gp}}$
lies in $P_0^{\textup{gp}}$.
To this aim, it suffices to show that$(\{ 0\} \oplus (\NN^{d-r})^{\textup{gp}})\cap P^{\textup{gp}} \subset P_0^{\textup{gp}}$.
Let $C(P)\subset P^{\gp}\otimes_{\ZZ}\QQ$ and $C(P_0)\subset P^{\gp}\otimes_{\ZZ}\QQ$ be cones spanned by $P$ and $P_0$
respectively.
Then $C(P)\cap P^{\gp}=P$ (cf. Remark~\ref{monoid1}), and the cone $C(P_0)$ is a face of $C(P)$.
Indeed, identifying $P^{\gp}\otimes_{\ZZ}\QQ$ with $F^{\gp}\otimes_{\ZZ}\QQ$,
the cone $C(P)$ and $C(P_0)$ are generated by
irreducible elements of $F=\NN^d$ and $\{ 0\}\oplus \NN^{d-r}$ respectively.
For any $p\in C(P_0)\cap P^{gp}$ there exists a positive integer $n$
such that $n\cdot p$ lies in $P_0$. Taking account of the definition of $P_0$
and $C(P_0)\cap P^{\gp}\subset P$,
we have $p\in P_0$, and thus $C(P_0)\cap P^{\gp}=P_0$.
Since $C(P_0)$ is a cone
in $P^{\gp}\otimes_{\ZZ}\QQ$, we have
$P_0^{\gp}\otimes_{\ZZ}\QQ\cap P^{\gp}=P_0^{\gp}$. (We regard
$P_0^{\gp}\otimes_{\ZZ}\QQ$ as a subspace of $P^{\gp}\otimes_{\ZZ}\QQ$.)
This means that
$(\{ 0\}\oplus (\NN^{d-r})^{\gp}\otimes_{\ZZ}\QQ)\cap P^{\gp}=P_0^{\gp}$.
Thus we conclude that
$(\{ 0\} \oplus (\NN^{d-r})^{\textup{gp}})\cap P^{\textup{gp}}\subset P_0^{\textup{gp}}$.
Hence we complete the proof.
\QED

\begin{Lemma}
\label{subcorrespond}
Let $R$ be a ring.
Let $X=\Spec R [\sigma^{\vee}\cap M]$ be the toric variety
over $R$,
where $\sigma$ is a full-dimensional simplicial cone in $N_{\RR}$ \textup{(}$N=\ZZ^d$\textup{)}.
Let $\MMM_X$ denote the canonical log structure induced by
$\sigma^{\vee}\cap M\to R [\sigma^{\vee}\cap M]$
Let $\overline{\MMM}_{X,\bar{x}}$
be the stalk at a geometric point $\bar{x}\to X$, and 
let $\overline{\MMM}_{X,\bar{x}}\to F$ be the minimal free resolution.
Then there exists a natural bijective map
from  the set of irreducible elements of $F$ to the set of torus-invariant divisors on $X$
on which $\bar{x}$
lies.

In particular, if $\sigma^{\vee}\cap M\to H$ is the minimal free resolution,
then there exists a natural bijective map
from the set of irreducible elements of $H$ to $\sigma(1)$.
\end{Lemma}

\Proof
Without loss of generality we may suppose that
$\bar{x}$ lies on the subscheme
\[
\Spec R[\sigma^{\vee}\cap M]/(\sigma^{\vee}\cap M)\subset \Spec R[\sigma^{\vee}\cap M].
\]
Then we have $\overline{\MMM}_{X,\bar{x}}=\sigma^{\vee}\cap M$.
The set of torus-invariant divisors on which $\bar{x}$ lies
is $\{ V(\rho)\}_{\rho\in \sigma(1)}$, i.e.,
the set of rays of $\sigma$.
For each ray $\rho\in \sigma(1)$
the intersection $\rho^{\perp}\cap \sigma^{\vee}$
is a $(\dim \sigma -1)$-dimensional face.
Since $\sigma^{\vee}$ is simplicial,
there is a unique ray of $\sigma^{\vee}$ which does not lie in $\rho^{\perp}\cap \sigma^{\vee}$. We denote this ray by $\rho^{\star}$.
Then it gives rise to a bijective map $\sigma(1)\to \sigma^{\vee}(1)
;\ \rho\mapsto \rho^{\star}$.
By Proposition~\ref{MFR}, there is a
natural embedding $\sigma^{\vee}\cap M \hookrightarrow F \hookrightarrow \sigma^{\vee}$
and each irreducible element of $F$ lies on a unique ray of $\sigma^{\vee}$.
It gives a bijective map from the set of irreducible elements of $F$
to $\sigma^{\vee}(1)$.
Hence our assertion follows.
\QED

\begin{Definition}
\label{admFR}
Let $S$ be a scheme. Let $N=\ZZ^d$ be a lattice and $M:=\Hom_{\ZZ}(N,\ZZ)$ its dual.
Let $\Sigma$ be a fan in $N_{\RR}$
and $X_{\Sigma}$ the associated toric variety
over $S$.
Let $\nn:=\{ n_\rho\}_{\rho\in\Sigma(1)}$ be a set of positive integers indexed by $\Sigma(1)$.
A morphism of fine log schemes $(f,\phi):(Y,\MMM)\to
(X_{\Sigma},\MMM_{\Sigma})$ is called an {\it admissible
FR morphism of type $\nn$}
if for any geometric point $\bar{y}\to Y$
the homomorphism
$f^{-1}\overline{\MMM}_{P,f(\bar{y})}\to \overline{\MMM}_{\bar{y}}$
is isomorphic to the composite $f^{-1}\overline{\MMM}_{P,f(\bar{y})}\stackrel{i}{\to} F\stackrel{n}{\to} F$ where $i:f^{-1}\overline{\MMM}_{P,f(\bar{y})}\to F$ is the minimal free resolution and $n:F\to F$ is defined by
$e_{\rho}\mapsto n_{\rho}\cdot e_{\rho}$.
Here for a ray $\rho\in \Sigma(1)$ such that $f(\bar{y})\in V(\rho)$
we write $e_{\rho}$ for the irreducible element of $F$ corresponding to
$\rho$
(cf. Lemma~\ref{subcorrespond}).
\end{Definition}

\begin{Proposition}
\label{ex2}
Let $P$ be a simplicially toric sharp monoid.
Suppose that $P=\sigma^{\vee}\cap M$ where $\sigma\subset N_{\RR}$
is a full-dimensional simplicial cone.
Let $\nn:=\{ n_\rho\}_{\rho\in\sigma(1)}$ be a set of positive integers indexed by $\sigma(1)$.
Let $\iota:P\to F$ be an admissible free resolution
defined to be the composite
$P\to F\stackrel{n}{\to} F$ where $P\to F$ is the minimal free resolution
and $n:F\to F$ is defined by $e_{\rho}\mapsto n_{\rho}\cdot e_{\rho}$
for each ray $\rho\in \sigma(1)$.
Here $e_{\rho}$ is the irreducible element of $F$ corresponding
to $\rho$ (cf. Definition~\ref{admtype}, Lemma~\ref{subcorrespond}).
Let $R$ be a ring.
Let $(f,\phi):(\Spec R[F],\MMM_F)\to (\Spec R[P],\MMM_P)$
be the  morphism induced by $\iota$.
Then $(f,\phi)$ is an admissible FR morphism of type $\nn$.
\end{Proposition}

\Proof
Let us denote by $(g,\psi):(\Spec R[F],\MMM_F)\to (\Spec R[F],\MMM_F)$
the morphism induced by $n:F\to F$.
Notice that for any geometric point $\bar{x}\to \Spec R[F]$
the canonical map $F\to \overline{\MMM}_{F,\overline{x}}$ and
$F\to \overline{\MMM}_{F,g(\bar{x})}$
are of the form $F\cong \NN^{s}\oplus \NN^t\stackrel{\textup{pr}_1}{\to}\NN^s\cong\overline{\MMM}_{F,\bar{x}}$
and $F\cong \NN^{s}\oplus \NN^t\stackrel{\textup{pr}_1}{\to}\NN^s\cong\overline{\MMM}_{F,g(\bar{x})}$ respectively for some $s,t\in\ZZ_{\ge 0}$.
Therefore $\overline{\psi}_{\bar{x}}:\NN^s\cong g^{-1}\overline{\MMM}_{F,g(\bar{x})}\to \NN^s\cong\overline{\MMM}_{F,\bar{x}}$ is the homomorphism induced by $e_{\rho}\mapsto n_{\rho}\cdot e_{\rho}$.
Taking account of Proposition~\ref{ex1},
our assertion follows from the definition of $\iota:P\to F$.
\QED

\begin{Proposition}
\label{CC}
Let $R$ be a ring.
Let $P$ be a simplicially toric sharp monoid and
$\iota:P\to F=\NN^d$ an admissible free resolution of $P$.
Let
$(f,h):(S,\MMM)\to (X_P:=\Spec R[P], \MMM_{P})$ be a morphism
of fine log schemes.
Here $\MMM_P$ is the canonical log structure on
$\Spec R[P]$ induced by $P\to R[P]$. Let $c:P \to \MMM_P$ be the standard chart.
Let $\bar{s}\to S$ be a geometric point.
Consider the following commutative diagram
\[
 \xymatrix{
 P \ar[d]_{\bar{c}_{\bar{s}}} \ar[r]^{\iota} & F \ar[d]^{\alpha}   \\
 f^{-1}\overline{\MMM}_{P,\bar{s}} \ar[r]^{\bar{h}_{\bar{s}}} & \overline{\MMM}_{\bar{s}} \\
  }
\]
where $\bar{c}_{\bar{s}}$ is the map induced by the standard chart
such that $\alpha$ \'etale locally lifts to a chart.
Then there exists an fppf neighborhood $U$ of
$\bar{s}$ in which we have a chart $\varepsilon :F\to \MMM$ such that
the following diagram
\[
 \xymatrix{
 P \ar[d]_{c} \ar[r]^{\iota} & F \ar[d]^{\varepsilon}   \\
 f^*\MMM_{P} \ar[r]^{{h}} & \MMM \\
  }
\]
commutes and the composition $F \stackrel{\varepsilon}{\to} \MMM \to \overline{\MMM}_{\bar{s}}$ is equal to $\alpha$.
If the order of the cokernel of $P^{\textup{gp}}\to F^{\gp}$
is invertible on $R$, then we can take an \'etale neighborhood $U$
of $\bar{s}$ with the above property.
\end{Proposition}

\Proof
Let $\gamma:=c_{\bar{s}}:P \to
f^{*}{\MMM}_{P,\bar{s}}$ be the chart induced by the standard chart $P\to \MMM_P$.
In order to show our assertion,
we shall prove that there exists a homomorphism 
$t:F \to \MMM_{\bar{s}}$, which is a lifting of $\alpha$,
such that $t\circ \iota=h_{\bar{s}}\circ \gamma$.
To this aim, consider the following diagram
\[
\xymatrix@R=3mm @C=12mm{
   & & \MMM_{\bar{s}}  \ar[dd]& &F \ar[dd]^{\xi} & \\
   &f^*\MMM_{P,\bar{s}}  \ar[dd] \ar[ur]^{h}& & P \ar[ll]^(0.3){\gamma} \ar[dd] \ar[ur]^(0.4){\iota}& & \\
   & &\MMM_{\bar{s}}^{\gp}   & &F^{\gp} \ar[r]_{\psi}^{\sim}  & \ZZ^d\\
    &f^*\MMM_{P,\bar{s}}^{\gp}  \ar[ru]^{h^{\gp}} & & P^{\gp} \ar[ll]^{\gamma^{\gp}} \ar[ru]^{\iota^{\gp}} \ar[r]_{\phi}^{\sim} &\ZZ^d \ar[ru]_{a} & \\
  }
\]
where vertical arrows are natural inclusions and $\phi$ and $\psi$
are isomorphisms chosen as follows.
By elementary algebra, we can take the isomorphisms $\phi$ and $\psi$  so that
$a:=\psi\circ \iota^{\gp}\circ \phi^{-1}$ is represented by the 
$(d\times d)$-matrix $(a_{ij})$
with $a_{ii}=:\lambda_i\in \ZZ_{\ge 1}$  for $0\le i\le d$,
and $a_{ij}=0$ for
$i\neq j$.
Let us construct a homomorphism $F^{\gp}\to \MMM_{\bar{s}}^{\gp}$
filling in the diagram.
Let $\{ e_i\}_{1\le i\le d}$ be the canonical basis of $\ZZ^d$
and put $m_i:=h^{\gp}\circ \gamma^{\gp}\circ \phi^{-1}(e_i)$.
Let $n_i$ be an element in $\MMM_{\bar{s}}^{\gp}$ such that
$q(n_i)=\alpha^{\gp}(\psi^{-1}(e_i))$ for $0\le i\le d$,
where $q:\MMM^{\gp}_{\bar{s}}\to \overline{\MMM}^{\gp}_{\bar{s}}$
is the natural projection. Then
there exists a unit element $u_i$ in $\OO_S^*$ such that
$ u_i+\lambda_i\cdot n_i=m_i$ in $\MMM_{\bar{s}}^{\gp}$ for $0\le i\le d$.
The algebra
$\OO':=\OO_{S,\bar{s}}[T_1,\ldots,T_d]/(T_i^{\lambda_i}-u_i)_{i=1}^d$ is 
a finite flat $\OO_{S,\bar{s}}$-algebra. If
the order of cokernel of $P^{\gp}\to F^{\gp}$ is 
invertible on $R$, $\OO'$ is an \'etale $\OO_{S,\bar{s}}$-algebra.
After the base change to $\OO'$, there exists a homomorphism 
$\eta:F^{\gp} \to {\MMM}_{\bar{s}}^{\gp}$
which is an extension of $h^{\gp}\circ \gamma^{\gp}$.
Since $\MMM_{\bar{s}}=\MMM_{\bar{s}}^{\gp}\times_{\overline{\MMM}_{\bar{s}}^{{\gp}}}\overline{\MMM}_{\bar{s}}$, the map
$t:F\to \MMM_{\bar{s}}^{\gp}\times_{\overline{\MMM}_{\bar{s}}^{{\gp}}}\overline{\MMM}_{\bar{s}}=\MMM_{\bar{s}}$ defined by $m \mapsto ({\eta}(\xi(m)),\alpha(m))$
is a homomorphism which makes
the diagram 
\[
 \xymatrix{
 P \ar[d]_{\gamma} \ar[r]^{\iota} & F \ar[d]^{t}   \\
 f^*\MMM_{P,\bar{s}} \ar[r]^{h_{\bar{s}}} & \MMM_{\bar{s}}, \\
  }
\]
commutative.
Since $P$ and $F$ are finitely generated,
this diagram extends to a chart in
some fppf neighborhood of $\bar{s}$.
The last assertion immediately follows.
\QED


\section{Moduli stack of admissible FR morphisms into a toroidal embedding}

\subsection{Moduli stack of admissible FR morphisms}
Let $(X,U)$ be a toroidal embedding over a field $k$.
Let us denote by $I$ the set of irreducible components
of $X-U$.

Consider a triple $(X,U,\nn)$, where $\nn=\{ n_{i}\in \ZZ_{\ge 1}\}_{i\in I}$.
We shall refer to $(X,U,\nn)$
as a {\it toroidal embedding $(X,U)$ of level
$\nn$}.
Let $\MMM_{X}$ be the canonical log structure of $(X,U)$.
The pair $(X,\mathcal{M}_X)$ is a fine saturated log scheme (cf. section 1.3).
If we further assume that $(X,U)$ is a {\it good} toroidal embedding (cf. section 1.2),
then
for any point $x$ on $X$, the stalk $\overline{\mathcal{M}}_{X,\bar{x}}$ is
a {\it simplicially toric sharp monoid}.

\begin{Proposition}
\label{correspond}
Let $\overline{\MMM}_{X,\bar{x}}$
be the stalk at a geometric point $\bar{x}\to X$.
Let $\overline{\MMM}_{X,\bar{x}}\to F$ be the minimal free resolution.
Then there exists a natural map
from  the set of irreducible elements of $F$ to the set of irreducible components of $X-U$
in which $\bar{x}$
lies.
\end{Proposition}
\Proof
Set $P=\overline{\MMM}_{X,\bar{x}}$.
Note first that there exists the natural correspondence between
irreducible elements of $F$ and height-one prime ideals of $P$ (see section 2).
Therefore it is enough to show that
there exists a natural map from the set of height-one prime ideals
of $P$ to the set of irreducible components of $X-U$ on which $\bar{x}$ lies.
Let $\{ \pp_1,\ldots,\pp_n\}$ be the set of height-one prime ideals of $P$.
Let $D=X-U$ be the reduced closed subscheme (each component
is pure 1-codimensional) and
$I_D$ be the ideal sheaf associated to $D$.
If $x$ denotes the image of $\bar{x}$, and $\OO_{X,x}$ denotes
its Zariski stalk, then the natural morphism $\Spec \OO_{X,\bar{x}}/I_D\OO_{X,\bar{x}}\to\Spec \OO_{X,x}/I_D\OO_{X,x}$ maps generic points to generic points.
On the other hand, the log structure $\MMM_{X}$ (cf. section 1)
has a chart $P\to \OO_{X,\bar{x}}$
and the support of $\overline{\MMM}_X$ is equal to $D$,
the underlying space of $\Spec \OO_{X,\bar{x}}/I_{D}\OO_{X,\bar{x}}$
is naturally equal to $\Spec \OO_{X,\bar{x}}/(\pp_1\OO_{X,\bar{x}}\cap\ldots\cap\pp_n\OO_{X,\bar{x}})$. (The closed subscheme
$\Spec k[P]/(\pp_1k[P]\cap\ldots \cap \pp_nk[P])$ has the same underlying space
as the complement $\Spec k[P]-\Spec k[P^{\gp}]$.)
Therefore it suffices to prove that
for any height-one prime ideal $\pp_i$, the closed subset
$\Spec \OO_{X,\bar{x}}/\pp_i\OO_{X,\bar{x}}$ is an irreducible
component of $\Spec \OO_{X,\bar{x}}/I_D\OO_{X,\bar{x}}$.
By \cite[Theorem 3.2 (1)]{log2}
there exists an isomorphism
$\hat{\OO}_{X,\bar{x}}\stackrel{\sim}{\to}k'[[P]][[T_1,\ldots,T_r]]$
and the composite $P\to \OO_{X,\bar{x}}\to \hat{\OO}_{X,\bar{x}}$
can be identified with the natural map $P\to k'[[P]][[T_1,\ldots,T_r]]$,
where $k'$ is the residue field of $\OO_{X,\bar{x}}$.
(Strictly speaking, \cite{log2} only treats the case of Zariski log structures,
but the same proof can apply to the case of \'etale log structures.)
Since $P-\pp_i$ is a submonoid and moreover it is a toric monoid,
thus $\hat{\OO}_{X,\bar{x}}/\pp_i\hat{\OO}_{X,\bar{x}}$ is
isomorphic to
the integral domain $k[[P-\pp_i]][[t_1,\ldots,t_r]]$.
Note that $\hat{\OO}_{X,\bar{x}}$ is a flat $\OO_{X,\bar{x}}$-algebra and the natural map $\OO_{X,\bar{x}}\to\hat{\OO}_{X,\bar{x}}$
is injective,
thus $\OO_{X,\bar{x}}/\pp_i\OO_{X,\bar{x}}\to \hat{\OO}_{X,\bar{x}}/\pp_i\hat{\OO}_{X,\bar{x}}$ is injective. Therefore $\OO_{X,\bar{x}}/\pp_i\OO_{X,\bar{x}}$
is an integral domain. Hence $\Spec \OO_{X,\bar{x}}/\pp_i\OO_{X,\bar{x}}$
is irreducible.
Thus we obtain the natural map as desired.
\QED

\begin{Definition}
\label{admdefinition}
Let $(X,U, \nn=\{ n_{i}\in \ZZ_{\ge 1}\}_{i\in I})$ be a good toroidal embedding of level $\nn$ over $k$.
(We denote by $I$ the set of irreducible components
of $X-U$.)
Let $\MMM_X$ be the canonical log structure on $X$.
An {\it admissible FR morphism} to 
$(X,\MMM_X,\nn)$
(or $(X,U,\nn)$)
is a morphism $(f,\phi):(S,\MMM_S)\to (X,\MMM_X)$
of fine log schemes such that for any geometric point
$\bar{s}\to S$ the homomorphism $\overline{\phi}:f^{-1}\overline{\MMM}_{X,f(\bar{s})}\to \overline{\MMM}_{S,\bar{s}}$
is isomorphic to \[
f^{-1}\overline{\MMM}_{X,f(\bar{s})}\stackrel{\iota}{\to}F\stackrel{n}{\to} F,
\]
where 
$\iota$ is the minimal free resolution and $n$ is defined by
$e\mapsto n_{i(e)}\cdot e$ where $e$ is an irreducible element of $F$
and $i(e)$ is an irreducible component
of $X-U$ to which $e$ corresponds via Proposition~\ref{correspond}.
(We shall call such a resolution an {\it admissible free resolution of type}
$\nn$ at $f(\bar{s})$.)
\end{Definition}

We define a category $\XXX_{(X,U, \nn)}$
as follows.
The objects are admissible FR morphisms to $(X,\MMM_X,\nn)$.
A morphism $\{(f,\phi):
(S,\MMM)\to (X,\MMM_X)\} \to \{(g,\psi):
(T,\NNN)\to (X,\MMM_X)\}$ in $\XXX_{(X,U,\nn)}$
is a morphism of $(X,\MMM_{X})$-log schemes
$(h,\alpha):(S,\MMM)\to (T,\NNN)$
such that $\alpha:h^*\NNN\to \MMM$ is an isomorphism.
By simply forgetting log structures, we have a functor
\[
\pi_{(X,U,\nn)}:\XXX_{(X,U,\nn)}\to (X\textup{-schemes}),
\]
which makes $\XXX_{(X,U,\nn)}$ a fibered category over the category of 
$X$-schemes.
This fibered category is a stack with respect to
\'etale topology because it is fibered in groupoid and
log structures are pairs of \'etale sheaves and their homomorphism
on \'etale site. Furthermore according to
the fppf descent theory for fine log
structures (\cite[Theorem A.1]{OL}),
$\XXX_{(X,U,\nn)}$
is a stack with respect to the fppf topology.

\begin{Theorem}
\label{Main1}
Let $(X,U,\nn=\{ n_{i}\in\ZZ_{\ge1}\}_{i\in I})$ be a good toroidal embedding
of level $\nn$
over a field $k$.
Then the stack $\XXX_{(X,U,\nn)}$ is
a smooth algebraic stack of finite type over $k$
with finite diagonal. 
The functor $\pi_{(X,U,\nn)}:\XXX_{(X,U,\nn)}\to X$
is a coarse moduli map, which induces an isomorphism
$\pi_{(X,U,\nn)}^{-1}(U)\stackrel{\sim}{\to}U$.
If we suppose that $(X,U)$ is tame (cf. section 1.2), and
$n_i$ is prime to the characteristic of $k$ for all $i\in I$,
then the stack
$\XXX_{(X,U,\nn)}$ is
a smooth Deligne-Mumford stack of finite type and separated
over $k$.

In the case of $n_{i}=1$ for all $i\in I$,
we have the followings.

\begin{enumerate}
\renewcommand{\labelenumi}{(\theenumi)}

\item The coarse moduli map induces an isomorphism
$\pi_{(X,U,\nn)}^{-1}(X_{\sm})\stackrel{\sim}\to X_{\sm}$ when
we denote by $X_{\sm}$ the smooth locus of $X$.

\item Suppose further that $(X,U)$ is tame.
If there exists another functor $f:\mathcal{X}\to X$
such that $\mathcal{X}$ is a smooth  separated Deligne-Mumford stack
and $f$ is a coarse moduli map, then there exists
a functor $\phi:\mathcal{X}\to\XXX_{(X,U,\nn)}$ such that the diagram
\[
 \xymatrix{
 \mathcal{X} \ar[r]^(0.3){\phi} \ar[dr]_{f} & \XXX_{(X,U,\nn)} \ar[d]^{\pi_{(X,U,\nn)}}   \\
  & X \\
  }
\]
commutes in the 2-categorical sense and such $\phi$ is unique up to a unique isomorphism.

\end{enumerate}

\end{Theorem}

\begin{Remark}
Moreover the stack $\XXX_{(X,U,\nn)}$ has the following nice properties,
which we will show later (because we need some preliminaries).

\begin{enumerate}
\renewcommand{\labelenumi}{(\theenumi)}

\item We will see that $\XXX_{(X,U,\nn)}$
is a {\it tame} algebraic stack in the sense of \cite{AOV}
(see Corollary~\ref{tamestack}).

\item We will prove that the complement $\XXX_{(X,U,\nn)}-U$
with reduced induced stack structure is a {\it normal crossing divisor}
on $\XXX_{(X,U,\nn)}$ (see section 3.5).

\end{enumerate}
\end{Remark}

\subsection{}Before the proof of Theorem~\ref{Main1},
we shall observe the case when
$(X,U)$ is an {\it affine simplicial toric variety} $\Spec R [\sigma^{\vee}\cap M]$
over a ring $R$,
where $\sigma$ is a full-dimensional simplicial cone in $N_{\RR}$ ($N=\ZZ^d$).
(For the application in section 4, we work over a general ring
rather than a field.)
Each ray $\rho\in \sigma$ defines the torus-invariant divisor $V(\rho)$.
Consider the torus embedding $X(\sigma,\nn):=(\Spec R [\sigma^{\vee}\cap M],\Spec
R[M], \nn)$ of level $\nn=\{ n_{\rho}\}_{\rho\in \sigma(1)}$.
Set $P:=\sigma^{\vee}\cap M$ (this is a simplicially toric sharp monoid).
Let $\iota:P\to F$ be the injective homomorphism
defined to be the composite $P\stackrel{i}{\to}F\stackrel{n}{\to}F$,
where $i$
is the
minimal free resolution
and $n:F\to F$ is defined by $e_{\rho}\mapsto n_{\rho}\cdot e_{\rho}$,
where $e_{\rho}$ denotes the irreducible element of $F$ corresponding
to a ray $\rho\in\sigma(1)$ (cf. Lemma~\ref{subcorrespond}).
The cokernel $F^{\gp}/\iota(P)^{\gp}$ is a finite abelian group
(and isomorphic to $F/\iota(P)$ by Lemma~\ref{monogp}).
We  view $F^{\gp}/\iota(P)^{\gp}$ as a finite commutative group scheme over $R$.
We denote by $G:=(F^{\gp}/\iota(P)^{\gp})^{D}$ its Cartier dual over $R$.
We define an action $m:\Spec R[F]\times_R G\to \Spec R[F]$ as follows.
Put $\pi:F\to F^{\gp}/\iota(P)^{\gp}$ the canonical map.
For each $R$-algebra $A$ and $A$-valued point $p:F\to A$ (a map of monoids)
of $\Spec R[F]$,
a $A$-valued point $g:(F^{\gp}/\iota(P)^{\gp})\to A$ (a map of monoids)
of $G$ sends $p$ to
the map $p^{g}:F\to A$ defined by the formula $p^g(f)= p(f)\cdot g(\pi(f))$
for $f \in F$.
We denote by $[\Spec R[F]/G]$ the stack-theoretic quotient
associated to the groupoid $m,\textup{pr}_1:\Spec R[F]\times_{R}G\rightrightarrows \Spec R[F]$ (cf. \cite[(10.13.1)]{LM}).
By \cite[section 3 and Theorem 3.1]{Con}, the coarse moduli space for this quotient stack is 
$\Spec R[F]^{G}$
where
\[
R[F]^{G}=\{ a \in R[F]|\ m^*(a)=\textup{pr}_1^*(a)\}\subset R[F].
\]
Moreover the natural morphism $R[P]\to R[F]^{G}$ is an isomorphism:

\begin{Claim}
\label{invariant}
The natural morphism $R[P]\to R[F]^{G}$ is an isomorphism.
In particular, the toric variety $\Spec R[P]$ is a coarse moduli space for
$[\Spec R[F]/G]$.
\end{Claim}

\Proof
Note first that $\Gamma (G)=R[F^{\gp}]/(f-1)_{f\in P^{\gp}\subset F^{\gp}}$
is a finite free $R$-module, and
\[
m^*:R[F]\to R[F]\otimes_{R}R[F^{\gp}]/(f-1)_{f\in P^{\gp}\subset F^{\gp}}
\]
maps $f\in F$ to $f\otimes f$.
Since $\operatorname{pr}^*_1$ maps $f$ to $f\otimes 1$, thus
$m^*(\Sigma_{f\in F}r_{f}\cdot f)=\operatorname{pr}^*_1(\Sigma_{f\in F}r_{f}\cdot f)$ ($r_f\in R$) if and only if
$r_{f}=0$ for $f$ with $f\notin P^{\gp}\subset F^{\gp}$.
Therefore it suffices to prove that $P^{\gp}\cap F=P$.
Clearly, we have $P^{\gp}\cap F\supset P$, thus we will show
$P^{\gp}\cap F\subset P$.
For any $f\in P^{\gp}\cap F$,
there exists a positive integer $n$ such that $n\cdot f\in P$.
Since $P$ is saturated, we have $f\in P$. Hence our claim follows.
\QED

Let $\MMM_{P}$ be the canonical log structure on the toric variety $\Spec R[P]$.

\begin{Proposition}
\label{localalg}
There exists a natural isomorphism
$\Phi:[\Spec R[F]/G]\to \XXX_{X(\sigma,\nn)}$ over $\Spec R[P]$.
The composite $\Spec R[F]\to [\Spec R[F]/G]\to \XXX_{X(\sigma,\nn)}$
corresponds to the admissible FR morphism
$(\chi,\epsilon):(\Spec R[F],\MMM_F)\to (\Spec R[P],\MMM_P)$
which is induced by $\iota:P\to F$.
\end{Proposition}
\Proof
Let $[\Spec R[F]/G] \to \Spec R[P]$ be the natural coarse moduli map.
According to \cite[Proposition 5.20 and Remark 5.21]{OL},
the stack $[\Spec R[F]/G]$ over $\Spec R[P]$ is isomorphic
to the stack $\mathcal{S}$ whose fiber over $f:S\to \Spec R[P]$
is the groupoid of triples $(\NNN, \eta,\gamma)$,
where $\NNN$ is a fine log structure on $S$,
$\eta :f^*\MMM_P\to\NNN$ is a morphism of log structures,
and $\gamma:F\to \overline{\NNN}$ is a morphism,
which \'etale locally lifts to a chart, such that the diagram
\[
 \xymatrix{
 P \ar[d]_{\bar{c}} \ar[r]^{\iota} & F \ar[d]^{\gamma}   \\
 f^{-1}\overline{\MMM}_{P} \ar[r] & \overline{\NNN} \\
  }
\]
commutes. Here we denote by $c$ the standard chart
$P\to\MMM_P$.
Notice that the last condition implies that $\eta$ is an admissible
FR morphism to $(\Spec R[P],\MMM_P,\nn)$.
Indeed, by Proposition~\ref{CC}
there exists fppf locally a chart $\gamma':F\to \NNN$
such that $\gamma'$ induces $\gamma$, and $\gamma'\circ \iota:P\to F\to \NNN$ is
equal to $P\to f^*\MMM_P\to \NNN$.
This means that there exists fppf locally on $S$ a strict morphism
$(S,\NNN)\to (\Spec R[F],\MMM_F)$ such that the composite
$(S,\NNN)\to (\Spec R[F],\MMM_F)\to (\Spec R[P],\MMM_P)$
is equal to $(f,\eta):(S,\NNN)\to (\Spec R[P],\MMM_P)$.
By Proposition~\ref{ex2} $(\Spec R[F],\MMM_F)\to (\Spec R[P],\MMM_P)$
is an admissible morphism to $(\Spec R[P],\MMM_P,\nn)$,
thus so is $(f,\eta):(S,\NNN)\to (\Spec R[P],\MMM_P)$.
Therefore there exists a
natural morphism
$\Phi:[\Spec R[F]/G] \to \XXX_{X(\sigma,\nn)}$
which forgets the additional data of the map $\gamma:F\to \overline{\NNN}$.
To show that $\Phi$ is essentially surjective,
it suffices to see that every object in $\XXX_{X(\sigma,\nn)}$ is fppf locally 
isomorphic to the image of $\Phi$.
Let
$(f,\phi):(S,\NNN)\to (\Spec R[P],\MMM_P)$
be an admissible FR morphism to $(\Spec R[P],\MMM_P,\nn)$.
By Proposition~\ref{ex2}, for any geometric point $\bar{s}\to S$ 
the map
$f^{-1}\overline{\MMM}_{P,f(\bar{s})}\to \overline{\NNN}_{\bar{s}}$
is isomorphic to $g^{-1}\overline{\MMM}_{P,g(\bar{x})}\to \overline{\MMM}_{F,\bar{x}}$,
where $g:\Spec R[F]\to \Spec R[P]$ is induced by $\iota:P\to F$,
and $\bar{x}$ is a geometric point on $\Spec R[F]$ which is lying over $f(\bar{s})$.
Since $F\to \overline{\MMM}_{F,\bar{x}}$
has the form
$F\cong \NN^r\oplus \NN^{d-r}\stackrel{\textup{pr}_1}{\to} \NN^{r}$,
thus
$P\to f^{-1}\overline{\MMM}_{P,f(\bar{s})}\to \overline{\NNN}_{\bar{s}}$
has the form $P\stackrel{\iota}\to F\cong \NN^r\oplus \NN^{d-r}\stackrel{\textup{pr}_1}{\to} \NN^{r}$,
thus the essentially surjectiveness follows from Proposition~\ref{CC}.
Finally, we will prove that $\Phi$ is fully faithful.
To this end, it suffices to show that
given two objects $(h_1:f^*\MMM_P\to \NNN_1,{\gamma}_1:F\to\overline{\NNN}_1)$ and
$(h_2:f^*\MMM_P\to \NNN_2, \gamma_2:F\to \overline{\NNN}_2)$ in $[\Spec R[F]/G](S)$,
any morphism of log structures $\xi:\NNN_1\to\NNN_2$, such that
$\xi\circ h_1=h_2$, has the property that $\bar{\xi}\circ \gamma_1=\gamma_2$.
It follows from the fact that $P$ is close to $F$ via $\iota$
and every stalk of $\overline{\NNN}_2$ is
of the form $\NN^r$ for some $r\in \NN$.
Indeed for any $f\in F$, there exists a positive integer
$n$ such that $n\cdot f\in P$. Since $\xi\circ h_1=h_2$,
we have $\bar{\xi}\circ \gamma_1(n\cdot f)=\gamma_2(n\cdot f)$.
Since every stalk of $\overline{\NNN}_2$ is
of the form $\NN^r$ for some $r\in \NN$, we conclude that
$\bar{\xi}\circ\gamma_1(f)=\gamma_2(f)$.

Finally, we will prove the last assertion. By Proposition~\ref{ex2},
$(\chi,\epsilon)$ is an admissible FR morphism of type $\nn$.
By \cite[Proposition 5.20]{OL}, the projection
$\Spec R[F]\to [\Spec R[F]/G]$ amounts to the triple
$(\MMM_F, \epsilon:\chi^* \MMM_P\to \MMM_F, F\to \overline{\MMM}_F)$
over $\chi:\Spec R[F]\to \Spec R[P]$.
Hence $\Phi((\MMM_F, \epsilon:\chi^* \MMM_P\to \MMM_F, F\to \overline{\MMM}_F))
=(\MMM_F, \epsilon:\chi^* \MMM_P\to \MMM_F)$ and thus our claim follows.
\QED

\begin{Proposition}
With the same notation as above, we have followings:
\label{additionals}
\begin{enumerate}
\renewcommand{\labelenumi}{(\theenumi)}

\item $\XXX_{X(\sigma,\nn)}$ has finite diagonal.

\item The natural morphism $\XXX_{X(\sigma,\nn)}\to \Spec R[P]$
is a coarse moduli map.

\item $\XXX_{X(\sigma,\nn)}$ is smooth
over $R$.

\end{enumerate}
\end{Proposition}

\Proof
We first prove (1).
The base change of the diagonal map $[\Spec R[F]/G]\to [\Spec R[F]/G]\times_{R}[\Spec R[F]/G]$ by the natural morphism $\Spec R[F]\times_{R}\Spec R[F]\to [\Spec R[F]/G]\times_{R}[\Spec R[F]/G]$ is isomorphic to
$\Spec R[F]\times_{R}G\to \Spec R[F]\times_{R}\Spec R[F]$, which maps $(x,g)$ to $(x,x^{g})$.
Since $\operatorname{pr}_1:\Spec R[F]\times_{R}G\to \Spec R[F]$ is proper
and $\operatorname{pr}_1:\Spec R[F]\times_{R}\Spec R[F]\to \Spec R[F]$
is separated, thus $\Spec R[F]\times_{R}G\to \Spec R[F]\times_{R}\Spec R[F]$
is proper. Clearly, it is also quasi-finite, so it is a finite morphism
and we conclude that $[\Spec R[F]/G]$ has finite diagonal.
Hence by Proposition~\ref{localalg} $\XXX_{X(\sigma,\nn)}$ has finite diagonal.
The assertion (2) follows from Claim~\ref{invariant}
and Proposition~\ref{localalg}.
Next we will prove (3).
By
Proposition \ref{localalg}, we have a finite flat cover
$\Spec R[F]\to \XXX_{X(\sigma,\nn)}$
from a smooth $R$-scheme $\Spec R[F]$, where $F$ is isomorphic to
$\NN^r$ for some $r\in\NN$.
Let $V\to \XXX_{X(\sigma,\nn)}$ be a smooth surjective morphism
from a $R$-scheme $V$.
Notice that the composite
$V\times_{\XXX_{X(\sigma,\nn)}}\Spec R[F] \stackrel{\textup{pr}_1}{\to} V\to \Spec R$
is smooth, and $V\times_{\XXX_{X(\sigma,\nn)}}\Spec R[F] \to V$ is a finitely presented flat surjective morphism.
Indeed the composite $V\times_{\XXX_{X(\sigma,\nn)}}\Spec R[F] \to \Spec R[F]\to\Spec R$ is smooth.
Thus by \cite[IV Proposition 17.7.7]{EGA},
we see that $V$ is smooth over $R$. Hence $\XXX_{X(\sigma,\nn)}$
is smooth over $R$.
\QED

\begin{Remark}
\label{fullgene}
Let $\sigma$ be a simplicial (not necessarily full-dimensional) cone in $N_{\RR}$.
Then there exists a splitting $N\cong N'\oplus N''$
such that $\sigma\cong \sigma'\oplus\{ 0\} \subset N'_{\RR}\oplus N''_{\RR}$
where $\sigma'$ is a full-dimensional cone in $N'_{\RR}$.
Thus $X_{\sigma}\cong X_{\sigma'}\times_{R}\Spec R[M'']$,
where $M''=\Hom_{\ZZ}(N'',\ZZ)$.
The log structure $\MMM_{\sigma}$ is isomorphic to the pullback
$\operatorname{pr}_1^*\MMM_{\sigma'}$ where $\operatorname{pr}_1: X_{\sigma'}\times_{R}\Spec R[M'']\to X_{\sigma}$.
Therefore $\XXX_{X(\sigma,\nn)}$ is an algebraic stack of finite type over $R$
for arbitrary simplicial cone $\sigma$. Also, Proposition~\ref{additionals}
holds for arbitrary simplicial cones.
\end{Remark}

\subsection{Proof of Theorem~\ref{Main1}}
We will return to the proof of Theorem~\ref{Main1}. Let $R=k$ be a field.
\smallskip

{\it Proof of algebraicity.}
We will prove that
$\XXX_{(X,U,\nn)}$ is
a smooth algebraic stack of finite type over a field $k$
with finite diagonal.
Clearly, one can assume that
$X(\sigma,\nn):=(X,U,\nn)=
(\Spec R [\sigma^{\vee}\cap M],\Spec
R[M], \nn)$ (with the same notation and hypothesis as in section 3.2),
where $\sigma$ is a full-dimensional simplicial cone.
Then by Proposition~\ref{localalg} and Proposition~\ref{additionals} (1) and (3), $\XXX_{X(\sigma,\nn)}$ is
a smooth algebraic stack of finite type over $R$ with finite diagonal
because $\XXX_{X(\sigma,\nn)}$
is isomorphic to $[\Spec\ R[F]/G]$.
Since the restriction $\MMM_X|_U$ is a trivial log structure,
thus we see that $\pi_{(X,U,\nn)}^{-1}(U)\to U$ is an {\it isomorphism}.

If we suppose that $(X,U)$ is tame and
$n_i$ is prime to the characteristic of $k$ for all $i\in I$,
then the order of $G=F^{\gp}/\iota(P)^{\gp}$ 
is prime to the characteristic of $k$. Therefore
the Cartier dual $(F^{\gp}/\iota(P)^{\gp})^{D}$
is a finite \'etale group scheme over $k$.
Thus $[\Spec k[F]/G]$ is a {\it Deligne-Mumford stack} over $k$
(cf. \cite[Remarque 10.13.2]{LM}).
Therefore in this case $\XXX_{(X,U,\nn)}$ is a Deligne-Mumford stack.

\vspace{1mm}

{\it Coarse moduli map for $\XXX_{(X,U,\nn)}$.}
Next we will prove that the natural map
$\XXX_{(X,U,\nn)}\to X$ is a coarse moduli map for $\XXX_{(X,U,\nn)}$.
By the above argument, we see that $\XXX_{(X,U,\nn)}\to X$
is a proper quasi-finite surjective morphism. Indeed, in the case when $X$ is a toric variety,
$\XXX_{(X,U,\nn)}\to X$ is a coarse moduli map by Proposition~\ref{additionals} (2). Moreover $\XXX_{(X,U,\nn)}$
is integral and $\XXX_{(X,U,\nn)}\to X$ is generically an isomorphism.
By \cite[Corollary 2.9 (ii)]{Hom},
we conclude that $\XXX_{(X,U,\nn)}\to X$ is a coarse moduli map.
\QED

Assume that $n_{i}=1$ for all $i\in I$.
Before the proof of (1) of Theorem~\ref{Main1},
we prove the following Lemma.

\begin{Lemma}
\label{sm}
Let $(X,U)$ be a toroidal embedding.
Let $\MMM_X$ be the canonical log structure of $(X,U)$.
For every geometric point $\bar{x}$ on $X_{\sm}$,
the monoid $\overline{\MMM}_{X,\bar{x}}$ is free.
\end{Lemma}

\Proof
Clearly, our claim is local with respect to \'etale topology.
Hence,
we may assume that for every closed point $x$ of $X$, there exist
an \'etale neighborhood $W\to X$, a simplicially toric sharp monoid $P$,
and a smooth morphism of $w:W\to \Spec k[P]$.
Note that $\MMM_X\mid_W$ is induced by the log structure on $\Spec k[P]$
defined by $P\to k[P]$.
Since $X$ is smooth, thus after shrinking $\Spec k[P]$ we may assume that
$\Spec k[P]$ is a smooth affine toric variety over $k$.
Therefore $\Spec k[P]$
has the form $\Spec k[\NN^r\oplus \ZZ^s]$ (cf. \cite[p28]{F})
and the log structure is induced by $\NN^r\to k[\NN^r\oplus \ZZ^s]$.
Hence our claim is clear.
\QED

 {\it Proof of (1) in Theorem~\ref{Main1}.}
To prove that $\pi_{(X,U,\nn)}^{-1}(X_{\sm})\stackrel{\sim}{\to}X_{\sm}$,
it suffices to show that $\pi_{(X,U,\nn)}^{-1}(X_{\sm})$ is an algebraic spaces
because $\pi_{(X,U,\nn)}$ is a coarse moduli map.
Thus let us
show that every MFR morphism
$(f,h):(S,\NNN)\to (X_{\sm},\MMM_X\mid_{X_{\sm}})$ has no non-trivial automorphism in $\XXX_{(X,U,\nn)}$.
By Lemma~\ref{sm} (1), each stalk of $f^{-1}\overline{\MMM}_X$ is a free monoid.
Hence $h:f^*\MMM_X\to\NNN$ is an isomorphism
because $(f,h)$ is an MFR morphism and thus $\bar{h}:f^{-1}\overline{\MMM}_X\to\overline{\NNN}$ is an isomorphism. Therefore
it does not have a non-trivial automorphism
and thus the map $\pi_{(X,U,\nn)}^{-1}(X_{\sm})\to X_{\sm}$
is an isomorphism.
\QED

{\it Proof of (2) in Theorem~\ref{Main1}.}
Let $p:Z\to\mathcal{X}$ be an \'etale cover by a separated
scheme $Z$ (note that $Z\times_{\XX}Z$ is also a separated scheme
since $\XX$ has finite diagonal and thus $Z\times_{\XX}Z\to Z\times_{k}Z$
is finite).
Let $X_{\textup{sing}}$ be the singular locus of $X$ and
put $V:=Z-p^{-1}(f^{-1}(X_{\textup{sing}}))$.
Note that since $X$ is a normal variety,
the codimension of $X_{\textup{sing}}$ is bigger than 1.
Thus the codimension of $p^{-1}(f^{-1}(X_{\textup{sing}}))$ is
bigger than 1.
Let $\operatorname{pr_1},\operatorname{pr_2}:Z\times_\mathcal{X}Z 
\rightrightarrows Z$
be natural projections.
By (1) in Theorem~\ref{Main1}, $f\circ p\mid_V:V\to X$ and $f\circ p\circ \operatorname{pr_i}\mid_{V\times_{\mathcal{X}}V}:V\times_{\mathcal{X}}V \to X \ (i=1,2)$
are uniquely lifted to morphisms into $\XXX_{(X,U,\nn)}$.
We abuse notation and denote by $f\circ p\mid_V$ and $f\circ p\circ \operatorname{pr_i}\mid_{V\times_{\mathcal{X}}V}$ ($i=1,2$)
lifted morphisms into $\XXX_{(X,U,\nn)}$.
Then by the purity lemma due to Abramovich and Vistoli (\cite[3.6.2]{A} \cite[Lemma 2.4.1]{AV}),
$f\circ p\mid_V$ and $f\circ p\circ \operatorname{pr_i}\mid_{V\times_{\mathcal{X}}V}$ ($i=1,2$)
are extended to $Z$ and $Z\times_{\mathcal{X}}Z$ respectively.
These extensions are unique up to a unique isomorphism.
Since $f\circ p\circ \operatorname{pr_1}=f\circ p\circ \operatorname{pr_2}$,
there exists a unique morphism $\phi:\mathcal{X} \to \XXX_{(X,U,\nn)}$
such that $\pi_{(X,U,\nn)}\circ\phi\cong f$.
\QED

\subsection{Stabilizer group schemes of points on $\XXX_{(X,U,\nn)}$}
\label{stab}
In this subsection, we calculate the stabilizer group schemes
(i.e. automorphism group schemes)
of points on the stack $\XXX_{(X,U,\nn)}$
for a good toroidal embedding $(X,U)$ of level $\nn=\{n_i\}_{i\in I}$.
For the definition of points on an algebraic stack,
we refer to \cite[Chapter 5]{LM}.
In this subsection, by a geometric point to an algebraic stack $\XX$
we mean
a morphism $\Spec K\to \XX$ with an algebraically closed field $K$.

Let us calculate the stabilizer groups of points on $\XXX_{(X,U,\nn)}$.
Let $\bar{x}:\Spec K\to \XXX_{(X,U,\nn)}$ be a point on $\XXX_{(X,U,\nn)}$
with an algebraically closed field $K$.
Note that
the point $\bar{x}$ can be naturally regarded as the point on $X$
via the coarse moduli map.
Suppose that $\bar{x}:\Spec K\to \XXX_{(X,U,\nn)}$ corresponds to
an admissible FR morphism $(\pi_{(X,U,\nn)}\circ \bar{x},h):(\Spec K,\MMM_K)\to
(X,\MMM_X,\nn)$.
Thus $\textup{Isom}(\bar{x},\bar{x})$
is the group scheme over $K$,
which represents the contravariant-functor
\[
G:(K\textup{-schemes})\to (\textup{groups}),
\]

$\{v:S\to \Spec K\}\mapsto \{ \textup{the group of isomorphisms
of the log structure}\ v^*\MMM_K$ which are compatible with
$v^*h:(\pi_{(X,U,\nn)}\circ \bar{x}\circ v)^*\MMM_{X}\to v^*\MMM_{K}\}$.
(We shall call $\textup{Isom}(\bar{x},\bar{x})$ the stabilizer
(or automorphism)
group scheme of point $\bar{x}$.)
Set $P=\overline{(\pi_{(X,U,\nn)}\circ \bar{x})^*\MMM}_{X}$ and
$F=\overline{\MMM}_K$.
Since $K$ is algebraically closed, 
there exist isomorphisms $(\pi_{(X,U,\nn)}\circ \bar{x})^*\MMM_{X}\cong
K^*\oplus P$ and $\MMM_K\cong K^*\oplus F$.
Therefore we have $G(\{S\to \Spec K\})=\Hom_{\textup{group}}(F^{\gp}/P^{\gp},\Gamma(S,\OO_S^*))$.
(Note that by Lemma~\ref{monogp}
the natural map $F/P\to F^{\gp}/P^{\gp}$ is an isomorphism.)
Hence we conclude that $G$ is the Cartier dual of $F^{\gp}/P^{\gp}$
over $K$.
Thus we have the following result:

\begin{Proposition}
\label{order}
Let $\bar{x}:\Spec K\to \XXX_{(X,U,\nn)}$
be a point with an algebraically closed field $K$.
Set $\bar{y}=\pi_{(X,U,\nn)}(\bar{x})$.
Then the stabilizer group scheme of $\bar{x}$ is 
isomorphic to the Cartier dual of $(F)^{\gp}/\overline{\MMM}_{X,\bar{y}}^{\gp}$
over $K$,
where $\overline{\MMM}_{X,\bar{y}}\to F$
is an admissible free resolution of type $\{ n_{i}\}$
at $\bar{y}$ (cf. Definition~\ref{admdefinition}).
\end{Proposition}

{\it Tame algebraic stacks.}
The recent paper \cite{AOV} introduced
the notion of tame algebraic stacks,
which is a natural generalization of that of tame Deligne-Mumford stacks.
We will show that the moduli stack $\XXX_{(X,U,\nn)}$ is a
tame algebraic stack.
Let us recall the definition of tame algebraic stacks.
Let $\XX$ be an algebraic stack over a base scheme $S$.
Suppose that the inertia stack $\XX\times_{\XX\times_{S}\XX}\XX$ is
finite over $\XX$.
Let $\pi:\XX\to X$ be a coarse moduli map for $\XX$
(Keel-Mori theorem implies the existence).
Let $\operatorname{QCoh}\XX$ (resp. $\operatorname{QCoh}X$)
denote the abelian category of quasi-coherent sheaves on $\XX$
(resp. $X$).
The algebraic stack $\XX$ is said to be {\it tame} if
the functor $\pi_*:\operatorname{QCoh}\XX\to\operatorname{QCoh}X$
is exact.
By \cite[Theorem 3.2]{AOV}, $\XX$ is tame if and only if
for any geometric point $\Spec K\to \XX$ with an algebraically
closed field $K$, its stabilizer group is a {\it linearly reductive} group
scheme over $K$. By Proposition~\ref{order},
for any geometric point $\Spec K\to \XXX_{(X,U,\nn)}$ with an algebraically
closed field $K$, its stabilizer group is {\it diagonalizable},
thus we obtain:

\begin{Corollary}
\label{tamestack}
The algebraic stack $\XXX_{(X,U,\nn)}$
is a tame algebraic stack.
\end{Corollary}

\subsection{Log structures on $\XXX_{(X,U,\nn)}$}
Let $(X,U,\nn=\{ n_i\}_{i\in I})$ be a good toroidal embedding
of level $\nn$
and $\XXX_{(X,U,\nn)}$ the associated stack.
Let us define a canonical (tautological) log structure
on $\XXX_{(X,U,\nn)}$.
Let us denote by $\textup{Lis-\'et}(\XXX_{(X,U,\nn)})$
the {\it lisse-\'etale site} of $\XXX_{(X,U,\nn)}$.
(Recall the definition of the lisse-\'etale site.
The underlying category of $\textup{Lis-\'et}(\XXX_{(X,U,\nn)})$
is the full subcategory of $\XXX_{(X,U,\nn)}$-schemes
whose objects are smooth $\XXX_{(X,U,\nn)}$-schemes.
A collection of morphisms
$\{f_i:S_i\to S \}_{i\in I}$ of smooth $\XXX_{(X,U,\nn)}$-schemes
is a covering family in $\textup{Lis-\'et}(\XXX_{(X,U,\nn)})$
if the morphism
$\sqcup_i\ f_i:\sqcup_i\ S_i\to S$
is \'etale surjective.)
We define a log structure $\MMM_{(X,U,\nn)}$
as follows. Let $s:S\to \XXX_{(X,U,\nn)}$
be a smooth $\XXX_{(X,U,\nn)}$-scheme, i.e., an object
in $\textup{Lis-\'et}(\XXX_{(X,U,\nn)})$.
This amounts exactly to an admissible
FR morphism $(f,\phi):(S,\MMM_S)\to (X,\MMM_{X},\nn)$
such that $f=\pi_{(X,U,\nn)}\circ s$.
By attaching to $s:S\to \XXX_{(X,U,\nn)}\in \textup{Lis-\'et}(\XXX_{(X,U,\nn)})$ the log structure $\MMM_{S}$, we define a fine and saturated log structure $\MMM_{(X,U,\nn)}$.
We shall refer this log structure as
the {\it canonical log structure} on $\XXX_{(X,U,\nn)}$.
(For the notion of log structures on algebraic stacks,
we refer to \cite[section 5]{OL}.)
Moreover by considering the homomorphisms $\phi:f^*\MMM_X\to \MMM_S$,
we have a natural morphism of log stacks
\[
(\pi_{(X,U,\nn)},\Phi):(\XXX_{(X,U,\nn)},
\MMM_{(X,U,\nn)})\to (X,\MMM_{X}).
\]

\begin{Proposition}
\label{log1}
Let $(X,U,\nn)$ be a good toroidal embedding with
a level $\nn=\{ n_i \}_{i\in I}$, where $I$ is the set of
irreducible elements of $X-U$.
Then we have the followings.
\begin{enumerate}
\renewcommand{\labelenumi}{(\theenumi)}

\item 
The subscheme $\DDD=\XXX_{(X,U,\nn)}-U$
with reduced closed subscheme structure is
a normal crossing divisor.
The log structure $\MMM_{(X,U,\nn)}$
is isomorphic to the log structure arising from
the divisor $\DDD=\XXX_{(X,U,\nn)}-U$.

\item Suppose further that
$(X,U)$ is a tame toroidal embedding (cf. section 1.2)
and $n_i$ is prime to the characteristic of the base field
for all $i$.
The morphism $(\pi_{(X,U,\nn)},\Phi):(\XXX_{(X,U,\nn)},
\MMM_{(X,U,\nn)})\to (X,\MMM_{X})$
is a Kummer log \'etale morphism.

\end{enumerate}
\end{Proposition}

We postpone the proof of this Proposition, and
it will be given in section 4.3
because it follows from the case when $(X,U)$ is a toric variety.

\begin{Remark}
One may regard $\XXX_{(X,U,\nn)}$ as a sort of ``stacky toroidal
embedding" endowed with the log structure $\MMM_{(X,U,\nn)}$.
\end{Remark}


\section{Toric algebraic stacks}

We define {\it toric algebraic stacks}.

\subsection{Some combinatorics}
Let $N=\ZZ^d$ be a lattice and $M=\Hom_{\ZZ}(N,\ZZ)$
the dual lattice. Let $\langle\bullet,\bullet\rangle :M\times N\to \ZZ$ be
the dual pairing.

$\bullet$ A pair $(\Sigma, \nn=\{n_\rho\}_{\rho\in\Sigma(1)})$
is called a {\it simplicial fan with a level structure $\nn$}
if $\Sigma$ is simplicial fan in $N_{\RR}$ and 
$\nn=\{n_\rho\}_{\rho\in\Sigma(1)}$ is the set of
positive integers indexed by the set of rays $\Sigma(1)$.

$\bullet$ A pair $(\Sigma,\Sigma^0)$ is called a {\it stacky fan}
if $\Sigma$ is simplicial fan in $N_{\RR}$ and
$\Sigma^0$ is a subset of $|\Sigma|\cap N$
such that
for any cone $\sigma$ in $\Sigma$ the restriction $\sigma\cap \Sigma^0$ is
a submonoid of $\sigma\cap N$ which has
the following properties:
(i) $\sigma\cap \Sigma^0$ is isomorphic to $\NN^r$ where $r=\dim \sigma$,
(ii) $\sigma\cap \Sigma^0$ is close to $\sigma\cap N$.
Put another way.
If $(\Sigma,\Sigma^0)$ is a stacky fan and $\rho$ is a ray of $\Sigma$,
then there exists the first point $w_{\rho}$ of
$ \rho\cap\Sigma^0$.
Since $\sigma\in \Sigma$ is simplicial and $\NN^{\dim \sigma}\cong\sigma\cap\Sigma^0$ is a submonoid that is close to $\sigma\cap N$,
thus $\sigma\cap\Sigma^0$ is the free monoid $\oplus_{\rho\in \sigma(1)}\NN\cdot w_{\rho}( \subset \sigma)$. (Each irreducible element of the monoid
$\sigma\cap\Sigma^0$ lies on a unique ray of $\sigma$.)
Therefore the data $\Sigma^0$ is determined by
the set of points $\{ w_{\rho}\}_{\rho\in\Sigma(1)}$.
Let us denote by $v_{\rho}$ the first lattice point on a ray
$\rho\in \Sigma(1)$.
For a ray $\rho\in \Sigma(1)$,
if $n_{\rho}\cdot v_{\rho}$ ($n_{\rho}\in \ZZ_{\ge 0}$)
is the first point $w_{\rho}$ of $\Sigma^0\cap \rho$,
then we shall call $n_{\rho}$ the {\it level of $\Sigma^0$ on $\rho$}.

$\bullet$ Let $(\Sigma, \nn=\{n_\rho\}_{\rho\in\Sigma(1)})$
be a simplicial fan with a level structure $\nn=\{n_{\rho}\}$.
The {\it free-net} of $\Sigma$ associated to level $\nn$
is the subset $\Sigma_{\nn}^{0}\subset |\Sigma| \cap N$
such that for each cone $\sigma\in \Sigma$
the set $\sigma\cap \Sigma_{\nn}^0$ is the free submonoid
generated by $\{n_{\rho}\cdot v_{\rho}\}_{\rho\in \sigma(1)}$,
where $v_{\rho}$ denotes the first lattice point on a ray $\rho\in \Sigma(1)$.
(This free-net may be regarded as the {\it geometric realization}
of the level $\nn$.)
The {\it canonical free-net} $\Sigma_{\textup{can}}^0$ of $\Sigma$
is the free-net associated to the level $\{ n_{\rho}=1\}_{\rho\in\Sigma(1)}$.

Note that for any stacky fan $(\Sigma,\Sigma^0)$
there exists a unique level $\nn=\{n_{\rho}\}_{\rho\in\Sigma(1)}$
such that $(\Sigma,\Sigma^0)=(\Sigma,\Sigma_{\nn}^0)$.

$\bullet$ Let $S$ be a scheme (or a ringed space). A stacky fan
$(\Sigma,\Sigma_{\nn}^0)$ is called {\it tame over} $S$
if for any cone $\sigma\in \Sigma$ the multiplicity $\mult (\sigma)$
is invertible on $S$, and for any $\rho\in \Sigma(1)$
the level $n_{\rho}$ is invertible on $S$.
For a stacky cone $(\sigma,\sigma_{\nn}^0)$,
we define the {\it multiplicity}, denoted by
$\mult (\sigma,\sigma_{\nn}^0)$, of $(\sigma,\sigma_{\nn}^0)$
to be $\mult(\sigma)\cdot \Pi_{\rho\in\sigma(1)}n_{\rho}$.
A stacky fan $(\Sigma,\Sigma_{\nn}^0)$
is tame over $S$ if and only if
$\mult (\sigma,\sigma_{\nn}^0)$ is invertible on $S$
for any cone $\sigma\in\Sigma$.

$\bullet$ A stacky fan
$(\Sigma,\Sigma_{\nn}^0)$ is called {\it complete} if
$\Sigma$ is a finite and complete fan, i.e., 
$\Sigma$ is a finite set and the support $|\Sigma|$ is
the whole space $N_{\RR}$.

\begin{Remark}
The notion of stacky fans was first introduced in \cite{BCS}.
\end{Remark}

\subsection{Toric algebraic stacks}

Fix a base scheme $S$.

\begin{Definition}
\label{toricstack}
Let $(\Sigma\subset N_{\RR},\nn=\{n_{\rho}\}_{\rho\in \Sigma(1)})$
be a simplicial fan with a level structure $\nn$.
Let $(\Sigma,\Sigma_{\nn}^0)$ be the associated stacky fan.
Define a fibered category
\[
\XXX_{(\Sigma,\Sigma_{\nn}^0)}\longrightarrow (S\textup{-schemes})
\]
as follows.
The objects over a $S$-scheme $X$ are triples
\[
(\pi:\SSS \to \OO_X,\alpha:\MMM\to \OO_X,\eta:\SSS\to \MMM)
\]
such that:

\begin{enumerate}
\renewcommand{\labelenumi}{(\theenumi)}

\item $\SSS$ is an \'etale sheaf of sub-monoids of the constant
sheaf $M$ on $X$ determined by $M=\Hom_{\ZZ} (N,\ZZ)$ such
that for every geometric point $\bar{x}\to X$, $\SSS_{x}\cong \SSS_{\bar{x}}$.
Here $x\in X$ is the image of $\bar{x}$
and $\SSS_{x}$ (resp. $\SSS_{\bar{x}}$) denotes the Zariski
(resp. \'etale) stalk. (The condition $\SSS_{x}\cong \SSS_{\bar{x}}$
for every $\bar{x}$ means that the \'etale sheaf $\SSS$ is 
arising from the Zariski sheaf $\SSS|_{X_{\operatorname{Zar}}}$.)

\item $\pi:\SSS\to \OO_X$ is a map of monoids where $\OO_X$ is a monoid
under multiplication.

\item For $s \in \SSS, \pi(s)$ is invertible if and only if
$s$ is invertible.

\item For each point $x \in X$, there exists some (and a unique)
$\sigma\in \Sigma$
such that $\SSS_{x}= \sigma^{\vee}\cap M$.

\item $\alpha:\MMM\to\OO_X$ is a fine log structure on $X$.

\item $\eta:\SSS \to \MMM$ is a homomorphism of sheaves of monoids
such that $\pi=\alpha\circ \eta$,
and for each geometric point $\bar{x}\to X$, the homomorphism $\bar{\eta}:\overline{\SSS}_{\bar{x}}
=(\SSS/(\textup{invertible\ elements}))_{\bar{x}}\to
\overline{\MMM}_{\bar{x}}$ is isomorphic to the composite
\[
\overline{\SSS}_{\bar{x}} \stackrel{r}{\hookrightarrow} F\stackrel{t}{\hookrightarrow} F,
\]
where $r$ is the minimal free resolution of $\overline{\SSS}_{\bar{x}}$
and $t$ is defined by
$e_{\rho} \mapsto n_{\rho}\cdot e_{\rho}$
where $e_{\rho}$ denotes the irreducible element of $F$
corresponding to a ray $\rho\in\Sigma(1)$ (see Lemma~\ref{defsub})
and  $n_{\rho}$ is the level of $\Sigma^0$ on $\rho$.
\end{enumerate}

A set of morphisms
\[
(\pi:\SSS \to \OO_X,\alpha:\MMM\to \OO_X,\eta:\SSS\to \MMM)\to(\pi':\SSS' \to \OO_X,\alpha':\MMM'\to \OO_X,\eta':\SSS'\to \MMM')
\]
over $X$
is the set of
isomorphisms of log structures $\phi:\MMM\to\MMM'$
such that $\phi\circ \eta=\eta':\SSS=\SSS'\to \MMM'$
if $(\SSS,\pi)=(\SSS',\pi')$
and
is the empty set if $(\SSS,\pi)\ne(\SSS',\pi')$.
With the natural notion of pullbacks, $\XXX_{(\Sigma,\Sigma_{\nn}^0)}$ is 
a fibered category.
According to \cite[Theorem on page 10]{AMRT}, for any $S$-scheme $X$,
there exists an isomorphism
\[
\Hom_{S\textup{-schemes}}(X,X_{\Sigma})\cong \{\textup{ all pairs }(\SSS,\pi)\textup{ on }X\textup{ satisfying}\ 
(1),(2),(3),(4) \},
\]
which commutes with pullbacks.
Here $X_{\Sigma}$ is the toric variety associated to $\Sigma$
over $S$.
Therefore there exists a natural functor 
\[
\pi_{(\Sigma,\Sigma_{\nn}^0)}:\XXX_{(\Sigma,\Sigma_{\nn}^0)} \longrightarrow X_{\Sigma}
\]
which simply forgets the data $\alpha:\MMM\to \OO_X$ and $\eta:\SSS\to \MMM$.
Moreover $\alpha:\MMM\to \OO_X$ and $\eta:\SSS\to \MMM$ are morphisms of the \'etale sheaves and thus $\XXX_{(\Sigma,\Sigma_{\nn}^0)}$ is a stack with respect to the
\'etale topology.

Objects of the form
$(\pi:M\to \OOO_X,\OOO_X^*\hookrightarrow \OOO_X,\pi:M\to \OOO_X^*)$
determine a substack of $\XXX_{(\Sigma,\Sigma_{\nn}^0)}$,
i.e., the natural inclusion
\[
i_{(\Sigma,\Sigma_{\nn}^0)}:T_{\Sigma}=
\Spec \OO_S[M]\hookrightarrow \XXX_{(\Sigma,\Sigma_{\nn}^0)}.
\]
We shall call $i_{(\Sigma,\Sigma_{\nn}^0)}$ the
{\it canonical torus embedding}.
This commutes with the torus-embedding
$i_{\Sigma}:T_{\Sigma}\hookrightarrow X_{\Sigma}$.
\end{Definition}

\begin{Lemma}
\label{defsub}
With the same notation as in Definition~\ref{toricstack},
let $e$ be an irreducible element
in $F$ and let $n$ be a positive integer such that
$n\cdot e\in r(\overline{\SSS}_{\bar{x}})$.
Let $m\in\SSS_{\bar{x}}$ be a lifting of $n\cdot e$.
Suppose that $\SSS_{\bar{x}}= \sigma^{\vee}\cap M$.
Then there exists a unique ray $\rho\in \sigma(1)$ such that
$\langle m,v_{\rho}\rangle>0$.
Here $v_{\rho}$
is the first lattice point of $\rho$,
and $\langle \bullet,\bullet\rangle$ is the dual pairing.
It does not depend on the choice of liftings.
Moreover this correspondence defines a natural injective map
\[
\{ \operatorname{Irreducible\ elements\ of } F \} \to \Sigma(1).
\]
\end{Lemma}

\Proof
Since the kernel of $\SSS_{\bar{x}}\to \overline{\SSS}_{\bar{x}}$
is $\sigma^{\perp}\cap M$, thus $\langle m,v_{\rho}\rangle$ does not depend upon the choice of liftings $m$.
Taking a splitting $N\cong N'\oplus N''$ such that
$\sigma\cong \sigma'\oplus \{ 0\}\subset N'_{\RR}\oplus N''_{\RR}$ where $\sigma'$ is a full-dimensional cone in
$N'_{\RR}$, we may and will assume that
$\sigma$ is a full-dimensional cone,
i.e., $\sigma^{\vee}\cap M$ is sharp.
By Proposition~\ref{AFR} (2),
there is a natural embedding $\sigma^{\vee}\cap M\hookrightarrow F\hookrightarrow \sigma^{\vee}$
such that each irreducible element
of $F$ lies on a unique ray of $\sigma^{\vee}$.
It gives rise to a bijective map from
the set of irreducible elements of $F$ to $\sigma^{\vee}(1)$. 
Since $\sigma$ and $\sigma^{\vee}$ are simplicial, we have a natural bijective map $\sigma^{\vee}(1)\to \sigma(1);
\ \rho\mapsto \rho^{\star}$,
where $\rho^{\star}$ is the unique ray which
does not lie in $\rho^{\perp}$.
Therefore the composite map from
the set of irreducible elements of $F$ to $\sigma(1)$
is a bijective map. Hence it follows our claim.
\QED

\begin{Remark}
\label{classic}
\begin{enumerate}
\renewcommand{\labelenumi}{(\theenumi)}

\item In what follows, we refer to
a homomorphism $\pi:\SSS\to \OO_X$ with properties
(1), (2), (3), (4) in Definition~\ref{toricstack}
as a {\it skeleton}. If $\pi:\SSS\to \OO_X$
corresponds to $X\to X_{\Sigma}$, then $\pi:\SSS\to \OO_X$
is called the {\it skeleton for} $X\to X_{\Sigma}$.

\item Let $(\Sigma,\Sigma_{\textup{can}}^0)$ be
a stacky fan such that $\Sigma$ is non-singular and
$\Sigma_{\textup{can}}^0$ is the canonical free net.
Then $\XXX_{(\Sigma,\Sigma_{\textup{can}}^0)}$
is the toric variety $X_{\Sigma}$ over $S$.
Indeed, for any object
$(\pi:\SSS \to \OO_X,\alpha:\MMM\to \OO_X,\eta:\SSS\to \MMM)$
in $\XXX_{(\Sigma,\Sigma_{\textup{can}}^0)}$
and any point $x\in X$,
the monoid $\overline{\SSS}_{\bar{x}}$
has the form $\NN^r$ for some $r\in \NN$.
Thus in this case, (5) and (6) in Definition~\ref{toricstack} are vacant.

\item Toric algebraic stack can be constructed over $\ZZ$
and pull back from there to any other scheme.
Therefore, for the proof of Proposition~\ref{toric1} and Theorem~\ref{maintoric} (1),
(2), (3), we may assume that the base scheme is $\Spec \ZZ$.

\end{enumerate}
\end{Remark}

{\it Torus Action functor.}
The torus action functor
\[
a:\XXX_{(\Sigma,\Sigma_{\nn}^0)} \times T_{\Sigma} \longrightarrow \XXX_{(\Sigma,\Sigma_{\nn}^0)}
\]
is defined as follows.
Let $\xi=(\pi:\SSS \to \OO_X,\alpha:\MMM\to \OO_X,\eta:\SSS\to \MMM)$
be an object in $\XXX_{(\Sigma,\Sigma_{\nn}^0)}$.
Let $\phi:M \to \OO_X$ be a map of monoids
from a constant sheaf $M$ on $X$ to $\OO_X$, i.e., an $X$-valued point
of $T_{\Sigma}:=\Spec \OO_S[M]$. Here $\OO_X$ is viewed as a sheaf of monoids
under multiplication.
We define $a(\xi,\phi)$ by 
$(\phi\cdot\pi:\SSS \to \OO_X,\alpha:\MMM\to \OO_X,\phi\cdot\eta:\SSS\to \MMM)$,
where $\phi\cdot\pi(s):=\phi(s)\cdot\pi(s)$ and
$\phi\cdot\eta(s):=\phi(s)\cdot\eta(s)$.
(Note that $\mathcal{S}$ is a subsheaf of the constant sheaf of $M$.)
Let $h:\MMM_1\to \MMM_2$ be a morphism in $\XXX_{(\Sigma,\Sigma_{\nn}^0)}\times T_{\Sigma}$
from $(\xi_1,\phi)$ to $(\xi_2,\phi)$, where
$\xi_i=(\pi:\SSS \to \OO_X,\alpha:\MMM_i\to \OO_X,\eta_i:\SSS\to \MMM_i)$
for $i=1,2$,
and $\phi:M \to \OO_X$ is an $X$-valued point
of $T_{\Sigma}$.
We define $a(h)$ to be $h$.
It gives rise to the functor $a:\XXX_{(\Sigma,\Sigma_{\nn}^0)} \times T_{\Sigma} \longrightarrow \XXX_{(\Sigma,\Sigma_{\nn}^0)}$ over $(S\textup{-schemes})$,
which makes the following diagrams
\[
\begin{CD}
\XXX_{(\Sigma,\Sigma_{\nn}^0)}\times T_{\Sigma}\times T_{\Sigma} @>{\textup{Id}_{\XXX_{(\Sigma,\Sigma_{\nn}^0)}}\times m}>> \XXX_{(\Sigma,\Sigma_{\nn}^0)}\times T_{\Sigma} \\
@VV{a\times \textup{Id}_{T_{\Sigma}}}V @VV{a}V \\
\XXX_{(\Sigma,\Sigma_{\nn}^0)}\times T_{\Sigma} @>{a}>> \XXX_{(\Sigma,\Sigma_{\nn}^0)}, \\
\end{CD}\qquad\qquad
\begin{CD}
\XXX_{(\Sigma,\Sigma_{\nn}^0)}\times T_{\Sigma} @>{a}>> \XXX_{(\Sigma,\Sigma_{\nn}^0)} \\
@A{\textup{Id}_{\XXX_{(\Sigma,\Sigma_{\nn}^0)}}\times e}AA @AA{\textup{Id}_{\XXX_{(\Sigma,\Sigma_{\nn}^0)}}}A \\
\XXX_{(\Sigma,\Sigma_{\nn}^0)} @= \XXX_{(\Sigma,\Sigma_{\nn}^0)}.
\end{CD}
\]
commutes in the strict sense, i.e., $a\circ (a\times \textup{Id}_{T_{\Sigma}})=a\circ
(\textup{Id}_{\XXX_{(\Sigma,\Sigma_{\nn}^0)}}\times m)$
and $a\circ (e\times \textup{Id}_{\XXX_{(\Sigma,\Sigma_{\nn}^0)}})
=\textup{Id}_{\XXX_{(\Sigma,\Sigma_{\nn}^0)}}$,
where $m:T_{\Sigma}\times T_{\Sigma}\to T_{\Sigma}$ is the natural
action and $e:S\to T_{\Sigma}$ is the unit section.
Thus the functor $a$ defines an action of $T_{\Sigma}$ on $\XXX_{(\Sigma,
\Sigma_{\nn}^0)}$ which extends the action of $T_{\Sigma}$ on itself
to the whole stack $\XXX_{(\Sigma,
\Sigma_{\nn}^0)}$
Here the notion of group actions on stacks is taken in the
sense of \cite[Definition 1.3]{R}.
This action makes the coarse moduli map $\pi_{(\Sigma,\Sigma_{\nn}^0)}:\XXX_{(\Sigma,\Sigma_{\nn}^0)}\to X_{\Sigma}$ torus-equivariant.

\begin{Proposition}
\label{toric1}
Let $(\Sigma,\nn=\{ n_{\rho}\}_{\rho\in\Sigma(1)})$
be a simplicial fan $\Sigma$ with level $\nn$.
Then there exists a canonical isomorphism between
the stack $\XXX_{(\Sigma,\Sigma_{\nn}^0)}$ and
the moduli stack $\XXX_{(X_{\Sigma},T_{\Sigma},\nn)}$ of 
admissible FR morphisms to $X_{\Sigma}$ of type $\nn$
\textup{(}cf. Definition~\ref{admFR} and section 3\textup{)}
over $X_{\Sigma}$.
Here $X_{\Sigma}$ is the toric variety over $S$.
\end{Proposition}

\Proof
We will explicitly
construct a functor $F:\XXX_{(X_{\Sigma},\nn)} \to \XXX_{(\Sigma,\Sigma_{\nn}^0)}$.
To this aim, consider the skeleton $\pi_u:\UU\to \OO_{X_{\Sigma}}$
for the identity element in $\Hom (X_{\Sigma},X_{\Sigma})$, i.e., the tautological object (cf. \cite[Theorem on page 10]{AMRT}).
Observe that the log structure associated to $\pi_u:\UU\to \OO_{X_{\Sigma}}$
is isomorphic to the canonical log structure $\MMM_{\Sigma}$.
Indeed let us
recall the construction of $\pi_u:\UU\to \OO_{X_{\Sigma}}$
(cf. \cite[page 11]{AMRT}).
By \cite[page 11]{AMRT}, we have
$\UU=\{$ union of subsheaves $(\sigma^{\vee}\cap M)_{X_{\sigma}}$ of the constant sheaf $M$ on $X_{\Sigma}$ of all $\sigma\in \Sigma \}$.
For a cone $\sigma\in \Sigma$, we have
$\UU(X_{\sigma})=\sigma^{\vee}\cap M$,
and $\UU(X_{\Sigma})\to \OO_{X_{\Sigma}}(X_{\sigma})$
is the natural map $\sigma^{\vee}\cap M\to \OO_S[\sigma^{\vee}\cap M]$.
Furthermore it is easy to see that if $V\subset X_{\sigma}$, then any element $m\in\UU(V)$
has the form $m_1+m_2$ where
$m_1\in \UU(X_{\sigma})$ and $m_2$ is an invertible element
of $\UU(V)$.
(If $\sigma\succ \tau$, then any element $m\in \tau^{\vee}\cap M$
has the form $m_1+m_2$ where
$m_1\in \sigma^{\vee}\cap M$ and $m_2$ is an invertible element
of $\tau^{\vee}\cap M$.)
Therefore $\pi_u:\UU\to \OO_{X_{\Sigma}}$ 
induces a homomorphism $\zeta:\UU \to \MMM_{\Sigma}$
which makes $\MMM_{\Sigma}$ the log structure associated to
$\UU\to \MMM_{\Sigma}$.
Let $(f,h):(X,\alpha:\MMM\to \OO_X) \to (X_{\Sigma},\MMM_{\Sigma})$
be an object in $\XXX_{(X_\Sigma,T_{\Sigma},\nn)}$ over $f:X\to X_{\Sigma}$.
Define $F((f,h):(X,\alpha:\MMM \to \OO_X) \to (X_{\Sigma},\MMM_{\Sigma}))$
to be
\[
(f^{-1}\pi_u :f^{-1}\UU \to \OO_X, \alpha:\MMM\to \OO_X,
h\circ (f^{-1}\zeta) :f^{-1}\UU \to \MMM).
\]
Then $(f^{-1}\pi_u :f^{-1}\UU \to \OO_X, \alpha:\MMM\to \OO_X,
h\circ (f^{-1}\zeta) :f^{-1}\UU \to \MMM)$ is an object in $\XXX_{(\Sigma,\Sigma_{\nn}^0)}$
because $(f,h)$ is an admissible FR morphism of type $\nn$
and the log structure $\MMM_{\Sigma}$ is arising from
the prelog structure $\pi_{u}:\UU\to \OO_{X_{\Sigma}}$.
Clearly, a morphism of log structure $\phi:\MMM\to\MMM$ commutes with
 $h\circ (f^{-1}\zeta) :f^{-1}\UU \to \MMM$ if and only if
 it commutes with $h:f^*\MMM_{\Sigma}\to \MMM$
 because the log structure the $f^*\MMM_{\Sigma}$
 is arising from the natural morphism
 $f^{-1}\UU \stackrel{f^{-1}\zeta}{\to} f^{-1}\MMM_{\Sigma} \to f^*\MMM_{\Sigma}$.
 Thus $F$ is fully faithful over $X\to X_{\Sigma}$.
 Finally, we will show that $F$ is essentially surjective over $X\to X_{\Sigma}$.
 Let $(f^{-1}\pi_u :f^{-1}\UU \to \OO_X, \alpha:\MMM\to \OO_X,
l :f^{-1}\UU \to \MMM)$ be an object of $\XXX_{(\Sigma,\Sigma_{\nn}^0)}$ over 
$f:X\to X_{\Sigma}$. (Note that every object in $\XXX_{(\Sigma,\Sigma_{\nn}^0)}$ is of this form.) 
  The homomorphism $l:f^{-1}\UU\to \MMM$ induces the homomorphism $l^a:f^*\MMM_{\Sigma} \to \MMM$.
Then $F((f,l^a):(X,\MMM)\to(X_{\Sigma},\MMM_{\Sigma}))$
 is isomorphic to $(f^{-1}\pi_u :f^{-1}\UU \to \OO_X, \alpha:\MMM\to \OO_X,
l :f^{-1}\UU \to \MMM)$.
Hence $F$ is an isomorphism.
\QED

\begin{Theorem}
\label{maintoric}
The stack $\XXX_{(\Sigma,\Sigma_{\nn}^0)}$
is a smooth tame algebraic stack 
that is locally of finite type over $S$
and has finite diagonal. Furthermore it satisfies the additional properties
such that:

\begin{enumerate}
\renewcommand{\labelenumi}{(\theenumi)}

\item The natural functor
$\pi_{(\Sigma,\Sigma_{\nn}^0)}:\XXX_{(\Sigma,\Sigma_{\nn}^0)} \longrightarrow X_{\Sigma}$ \textup{(}cf. Definition~\ref{toricstack}\textup{)}
is a coarse moduli map,

\item The canonical torus embedding $i_{(\Sigma,\Sigma_{\nn}^0)}:T_{\Sigma}=
\Spec \OO_S[M]\hookrightarrow \XXX_{(\Sigma,\Sigma_{\nn}^0)}$
is an open immersion identifying $T_{\Sigma}$
with a dense open substack of $\XXX_{(\Sigma,\Sigma_{\nn}^0)}$.

\item If $\Sigma$ is a finite fan, then $\XXX_{(\Sigma,\Sigma_{\nn}^0)}$
is of finite type over $S$.

\end{enumerate}

If $(\Sigma,\Sigma_{\nn}^0)$ is tame over the base scheme $S$,
then $\XXX_{(\Sigma,\Sigma_{\nn}^0)}$
is a Deligne-Mumford stack over $S$.
\textup{(}Moreover there is a criterion for Deligne-Mumfordness. see Corollary~\ref{Deligne-Mumfordness}.\textup{)}

If $\Sigma$ is a non-singular fan and $\Sigma_{\textup{can}}^0$
denotes the canonical free net, then $\XXX_{(\Sigma,\Sigma_{\textup{can}}^0)}$
is the toric variety $X_{\Sigma}$ over $S$.
\end{Theorem}

\Proof
By Proposition~\ref{toric1},
Proposition~\ref{localalg}, and Proposition~\ref{additionals} (see also
Remark~\ref{fullgene}),
$\XXX_{(\Sigma,\Sigma_{\nn}^0)}$ is a smooth algebraic stack
that is locally of finite type
over $S$ and has finite diagonal.
Moreover $\pi_{(\Sigma,\Sigma_{\nn}^0)}:\XXX_{(\Sigma,\Sigma_{\nn}^0)} \rightarrow X_{\Sigma}$ is a coarse moduli map by Proposition~\ref{additionals} (2), and thus by Keel-Mori theorem
(cf. \cite{KM}, see also \cite[Theorem 1.1]{Con})
$\pi_{(\Sigma,\Sigma_{\nn}^0)}$ is proper.
Hence if $\Sigma$ is a finite fan, then
$\XXX_{(\Sigma,\Sigma_{\nn}^0)}$ is of finite type over $S$.
To see that $\XXX_{(\Sigma,\Sigma_{\nn}^0)}$ is tame,
we may assume that the base scheme is a spectrum of a field
because
the tameness depends only on automorphism group schemes
of geometric points (cf. \cite[Theorem 3.2]{AOV}).
Thus the tameness follows from Proposition~\ref{toric1}
and Corollary~\ref{tamestack}.
Since the restriction $\MMM_{\Sigma}|_{\Spec \OO_S[M]}$ of
$\MMM_{\Sigma}$ to $T_{\Sigma}=\Spec \OO_S[M]$ is the trivial log structure,
thus by taking account of Proposition~\ref{toric1} the morphism
$\pi_{(\Sigma,\Sigma_{\nn}^0)}^{-1}(T_{\Sigma}) \to T_{\Sigma}$
is an isomorphism. Hence $i_{(\Sigma,\Sigma_{\nn}^0)}:T_{\Sigma}=
\Spec \OO_S[M]\to \XXX_{(\Sigma,\Sigma_{\nn}^0)}$
is
an open immersion.
If $(\Sigma,\Sigma_{\nn}^0)$ is tame over the base scheme $S$,
the same argument as in {\it Proof of algebraicity} of Theorem~\ref{Main1}
shows that $\XXX_{(\Sigma,\Sigma_{\nn}^0)}$ is a Deligne-Mumford stack.
The last claim follows from Remark~\ref{classic}.
\QED

\begin{Definition}
Let $(\Sigma,\Sigma_{\nn}^0)$ be a stacky fan.
We shall call the stack $\XXX_{(\Sigma,\Sigma_{\nn}^0)}$
the {\it toric algebraic stack} associated to $(\Sigma,\Sigma_{\nn}^0)$
(or $(\Sigma,\nn)$).

\end{Definition}

\begin{Remark}
\label{diag}
\begin{enumerate}
\renewcommand{\labelenumi}{(\theenumi)}

\item Let $(\textup{Smooth toric varieties over}\ S)$,
(resp. $(\textup{Simplicial toric varieties over}\ S)$)
denote the category of smooth (resp. simplicial) toric varieties
over $S$ whose morphisms are $S$-morphisms.
Let $(\textup{Toric algebraic stacks over}\ S)$
denote the 2-category whose objects are toric algebraic stacks
over $S$, and a 1-morphism is an $S$-morphism between objects,
and a 2-morphism is an isomorphism between 1-morphisms.
Given a 1-morphism $f:\XXX_{(\Sigma,\Sigma_{\nn}^0)}\to \XXX_{(\Delta,\Delta_{\nn'}^0)}$, by the universality of coarse moduli spaces,
there exists a unique morphism
$f_0:X_{\Sigma}\to X_{\Delta}$
such that $\pi_{(\Delta,\Delta_{\nn'}^0)}\circ f=f_0\circ \pi_{(\Sigma,\Sigma_{\nn}^0)}$.
By attaching $f_0$ to $f$, we obtain a functor
\[
c:(\textup{Toric algebraic stacks over}\ S)\to(\textup{Simplicial toric varieties over}\ S),\ \XXX_{(\Sigma,\Sigma_{\nn}^0)}\mapsto X_{\Sigma}.
\]
Therefore there is the following diagram of (2)-categories,
\[
\xymatrix@R=1mm @C=12mm{
  & (\textup{Toric algebraic stacks over}\ S)\ar[dd]^c \\
(\textup{Smooth toric varieties over}\ S)\ar[ur]^a \ar[dr]^b & \\
 &  (\textup{Simplicial toric varieties over}\ S) \\
}
\]
where $a$ and $b$ is fully faithful functors and
$c$ is an essentially surjective functor.

\item {\it Gluing pieces together}.
Let $(\Sigma,\Sigma_{\nn}^0)$ be a stacky fan
and $\sigma$ be a cone in $\Sigma$.
By definition there exists
a natural fully faithful morphism
$\XXX_{(\sigma,\sigma_{\nn}^0)}\hookrightarrow \XXX_{(\Sigma,\Sigma_{\nn}^0)}$
where $(\sigma,\sigma_{\nn}^0)$ is the restriction of
$(\Sigma,\Sigma_{\nn}^0)$ to $\sigma$.
The image of this functor is identified with the
open substack $\pi_{(\Sigma,\Sigma_{\nn}^0)}^{-1}(X_{\sigma})$,
where $X_{\sigma}( \subset X_{\Sigma})$ is the
affine toric variety associated to $\sigma$.
That is to say, Definition~\ref{toricstack} allows
one to have a natural gluing construction
$\cup_{\sigma\in \Sigma}\XXX_{(\sigma,\sigma_{\nn}^0)}=\XXX_{(\Sigma,\Sigma_{\nn}^0)}$.

\end{enumerate}
\end{Remark}

\begin{Proposition}
\label{toric2}
Let $(\sigma, \sigma_{\nn}^0)$
be a stacky fan such that $\sigma$ is a simplicial cone
in $N_{\RR}$ where $N=\ZZ^d$.
Here $N\cong\ZZ^{d}$.
Suppose that $\dim (\sigma)=r$.
Then the toric algebraic stack $\XXX_{(\sigma,\sigma_{\nn}^0)}$ has a finite
fppf morphism
\[
p:\mathbb{A}_{S}^{r}\times\mathbb{G}_{m,S}^{d-r}\to \XXX_{(\sigma,\sigma_{\nn}^0)}
\]
where $\mathbb{A}_S^r$ is an $r$-dimensional affine space
over $S$.
Furthermore $\XXX_{(\sigma,\sigma_{\nn}^0)}$ is isomorphic to the quotient stack
$[\mathbb{A}_S^r/G]\times \GG_m^{d-r}$ where $G$ is a finite
flat group scheme
over $S$. If $\sigma$ is a full-dimensional cone,
the quotient $[\mathbb{A}_S^r/G]$ coincides with
the quotient given in section 3.2 (cf. Proposition~\ref{localalg}).

If $(\sigma,\sigma_{\nn}^0)$ 
is tame over
the base scheme $S$,
then we can choose $p$ to be a finite \'etale cover and $G$ to be
a finite \'etale group scheme.
\end{Proposition}

\Proof
As in Remark~\ref{fullgene}
we choose a splitting
\begin{align*}
N &\cong N'\oplus N'' \\
(\sigma,\sigma_{\nn}^0) &\cong (\tau,\tau_{\nn}^0)\oplus \{ 0\} \\
X_{\sigma} &\cong X_{\tau}\times\mathbb{G}_{m,S}^{d-r}
\end{align*}
where $\sigma\cong \tau \subset N'_{\RR}$ is a full-dimensional cone.
Notice that $P:=\tau^{\vee}\cap M'$ is a simplicially toric sharp monoid.
Let $\iota:P\to F\cong\NN^r$ be the homomorphism of monoids
defined as the composite
$P\to F\stackrel{n}{\to} F$ where $P\to F$ is the minimal free resolution
and $n:F\to F$ is defined by $e_{\rho}\mapsto n_{\rho}\cdot e_{\rho}$
for each ray $\rho\in \sigma(1)$.
Here $e_{\rho}$ is the irreducible element of $F$ corresponding
to $\rho$ (cf. Lemma~\ref{subcorrespond}).
Then by Proposition~\ref{localalg} and Proposition~\ref{toric1},
the stack $\XXX_{(\tau,\tau_{\nn}^0)}$ is isomorphic to the quotient stack
$[\Spec \OO_S[F]/G]$ where $G$ is the Cartier dual $(F^{\gp}/\iota(P)^{\gp})^D$.
The group scheme $G$ is a finite flat over $S$.
We remark that $\XXX_{(\sigma,\sigma_{\nn}^0)}\cong\XXX_{(\tau,\tau_{\nn}^0)}\times\mathbb{G}_{m,S}^{d-r}$
by Proposition~\ref{toric1}, and thus our claim follows.
The last assertion is clear because in such case $G$ is a finite
\'etale group scheme.
\QED

\begin{Proposition}
\label{properness}
The toric algebraic stack $\XXX_{(\Sigma,\Sigma_{\nn}^0)}$ is proper over $S$
if and only if
$(\Sigma,\Sigma_{\nn}^0)$ is complete (cf. section 4.1).
\end{Proposition}

\Proof
If $\Sigma$ is a finite and complete fan, then since the coarse
moduli map $\pi_{(\Sigma,\Sigma_{\nn}^0)}$ is proper,
properness of $\XXX_{(\Sigma,\Sigma_{\nn}^0)}$ over $S$
follows from the fact that $X_{\Sigma}$ is proper over $S$ (cf. \cite[Proposition in page 39]{F} or \cite[Chapter IV, Theorem 2.5 (viii)]{CF}).
Conversely, suppose that $\XXX_{(\Sigma,\Sigma_{\nn}^0)}$ is proper over $S$.
We will show that $\Sigma$ is a finite fan and the support
$|\Sigma|$ is the whole space $N_{\RR}$.
By \cite[Lemma 2.7(ii)]{Hom}, we see that $X_{\Sigma}$
is
proper over $S$. By (\cite[Proposition in page 39]{F} or \cite[Chapter IV, Theorem 2.5 (viii)]{CF}),
it suffices to show only that $\Sigma$ is a finite fan (this
is well known, but we give the proof here because \cite{F} assumes
that all fans are finite).
Let $\Sigma_{\textup{max}}$ be the subset of $\Sigma$ consisting of
the maximal elements with respect to the face relation.
The set $\{ X_{\sigma}\}_{\sigma\in\Sigma_{\textup{max}}}$
is an open covering of $X_{\Sigma}$, but
$\{ X_{\sigma}\}_{\sigma\in\Sigma'}$ is not an open covering
for any proper subset $\Sigma'\subset \Sigma_{\textup{max}}$
because for any cone $\sigma$, set-theoretically $X_{\sigma}=\sqcup_{\sigma\succ \tau}Z_{\tau}$ (cf. section 1.1).
Since $X_{\Sigma}$ is quasi-compact, $\Sigma_{\textup{max}}$
is finite. Hence $\Sigma$ is a finite fan (note that $\Sigma$ is simplicial in our case).
\QED

\begin{Remark}
Here we will explain the relationship between toric Deligne-Mumford
stacks introduced in \cite{BCS} and toric algebraic stacks introduced
in this section.
We assume that the base scheme is an algebraically closed field
of characteristic zero (since the results of \cite{BCS} only hold in characteristic zero).
First, since both took different approaches, it is not clear that
if one begins with a given stacky fan, then
two associated stacks in the sense of \cite{BCS} and us are isomorphic to each other. But in the subsequent paper \cite{I2} we prove the {\it geometric characterization theorem} for toric algebraic stacks in the sense of us.
One can deduce from it that they are (non-canonically) isomorphic to each other.
(See \cite[Section 5]{I2}.)
Next, toric Deligne-Mumford stacks in the sense of \cite{BCS}
admits a finite abelian gerbe structure, while our toric algebraic stacks
do not.
However, such structures can be obtained form our toric algebraic stacks
by the following well-known technique. We here work over $\ZZ$.
Let $\LLL_{D}$ be an invertible sheaf on $\XXX_{(\Sigma,\Sigma_{\nn}^0)}$,
associated to a torus invariant divisor $D$ (see the next subsection).
Consider the triple
\[
(u:U\to \XXX_{(\Sigma,\Sigma_{\nn}^0)},\mathcal{M},\ \phi:u^*\LLL_D\cong \mathcal{M}^{\otimes n})
\]
where $\mathcal{M}$ is an invertible sheaf on $U$, and $\phi$ is an isomorphism
of sheaves.
Morphisms of triples are defined in the natural manner.
Then it forms an algebraic stack $\XXX_{(\Sigma,\Sigma_{\nn}^0)}(\LLL_D^{1/n})$
and there exists the natural forgetting functor 
$\XXX_{(\Sigma,\Sigma_{\nn}^0)}(\LLL_D^{1/n})\to \XXX_{(\Sigma,\Sigma_{\nn}^0)}$, which is a smooth morphism. The composition of this procedure, that is,
\[
\XXX_{(\Sigma,\Sigma_{\nn}^0)}(\LLL_{D_1}^{1/n_1})\times_{\XXX_{(\Sigma,\Sigma_{\nn}^0)}
}\cdots \times_{\XXX_{(\Sigma,\Sigma_{\nn}^0)}
}\XXX_{(\Sigma,\Sigma_{\nn}^0)}(\LLL_{D_r}^{1/n_r})\to \XXX_{(\Sigma,\Sigma_{\nn}^0)}
\]
yields a gerbe that appears on toric Deligne-Mumford stacks.
The stack $\XXX_{(\Sigma,\Sigma_{\nn}^0)}(\LLL_D^{1/n})$ associated to
$\LLL_D$ can be viewed as the fiber product
\[
\xymatrix{
\XXX_{(\Sigma,\Sigma_{\nn}^0)}(\LLL_D^{1/n}) \ar[r] \ar[d] & B\GG_m \ar[d] \\
\XXX_{(\Sigma,\Sigma_{\nn}^0)} \ar[r] & B\GG_m, \\
}
\]
where the lower horizontal arrow is associated to $\LLL_{D}$,
$B\GG_m$ is the classifying stack of $\GG_m$, and
$B\GG_m\to B\GG_m$ is associated to the homomorphism
$\GG_m\to \GG_m, a\mapsto a^n$.
This interpretation has the following direct generalization.
Let $N$ be a finitely generated abelian group
and $h:\ZZ\to N$ be a homomorphism of abelian
groups.
Then it gives rise to $G=\Spec \ZZ[N]\to \GG_m$ and $BG\to B\GG_m$.
Then we can define $\XXX_{(\Sigma,\Sigma_{\nn}^0)}(\LLL_D,h)$
to be the fiber product
\[
\xymatrix{
 & BG \ar[d] \\
\XXX_{(\Sigma,\Sigma_{\nn}^0)} \ar[r] & B\GG_m. \\
}
\]
Notice that $\XXX_{(\Sigma,\Sigma_{\nn}^0)}(\LLL_D,h)$
is smooth, but does not have the finite inertia stack, that is,
automorphism groups of objects in $\XXX_{(\Sigma,\Sigma_{\nn}^0)}(\LLL_D,h)$
can be positive dimensional.

\end{Remark}


\subsection{Torus-invariant cycles and log structures on toric 
algebraic stacks}

In this subsection we define torus-invariant cycles
and a canonical log structure
on toric algebraic stacks.
It turns out that
if a base scheme $S$ is regular the complement of a torus embedding
in a toric algebraic stack
is a {\it divisor with normal crossings} relative to $S$ (unlike simplicial
toric varieties),
and it fits nicely in the notions of torus-invariant divisors and 
canonical log structure.

{\it Torus-invariant cycles}. Let $(\Sigma \subset N_{\RR},\Sigma_{\nn}^0)$ be a stacky fan,
where $N=\ZZ^d$.
The torus-invariant cycle
$V(\sigma)\subset X_{\Sigma}$ associated to $\sigma$
is represented by the functor
\[
F_{V(\sigma)}:(S\textup{-schemes})\to (\textup{Sets}),
\]
\[
X\mapsto \{\textup{all pairs}\ (\SSS, \pi:\SSS\to \OO_X)\ \textup{of skeletons with property}
\ (*) \}
\]
where
$(*)=\{\textup{For any point}\ x\in X,\ \textup{there exists a cone}\ \tau\ \textup{such that}\ \sigma \prec \tau\ \textup{and}\ \SSS_{x}=\tau^{\vee}\cap M,\ 
\textup{and}\ \pi(s)=0\ \textup{if}\ s\in \tau^{\vee}\cap \sigma_o^{\vee}\cap M\subset \SSS_{x} \}$ (see section 1.1 for the definition of $\sigma_o^{\vee}$).
Indeed, to see this,
we may and will suppose that $X_{\Sigma}$ is an affine toric variety
$X_{\tau}=\Spec \OO_S[\tau^{\vee}\cap M]$ with $\sigma\prec \tau$.
Let $\pi:\UU\to \OO_{X_{\tau}}$ be the skeleton for the identity morphism $X_{\tau}\to X_{\tau}$.
As in the proof of Proposition~\ref{toric1},
$\UU=\{$ union of subsheaves $(\gamma^{\vee}\cap M)_{X_{\gamma}}$ of the constant sheaf $M$ on $X_{\tau}$ of all $\gamma\prec \tau \}$.
Let $\UU':=\{$ union of subsheaves $(\gamma^{\vee}\cap \sigma_o^{\vee}\cap M)_{X_{\gamma}}$ of the constant sheaf $M$ on $X_{\tau}$ of all $\gamma\prec \tau \}$.
Then $F_{V(\sigma)}$ is represented by 
$\Spec (\OO_S[\tau^{\vee}\cap M]/\pi(\UU'))$.
Since $\gamma^{\vee}\cap \sigma_o^{\vee}\cap M=\UU'(X_{\gamma})\subset\gamma^{\vee}\cap M=\UU(X_{\gamma})$ for $\gamma\prec \tau$,
it is easy to check that $\Spec (\OO_S[\tau^{\vee}\cap M]/\pi(\UU'))=\Spec (\OO_S[\tau^{\vee}\cap M]/\tau^{\vee} \cap\sigma_o^{\vee}\cap M)=V(\sigma)$.

Now let us define the torus-invariant cycle $\VVV(\sigma)$
associated to $\sigma\in \Sigma$.
Consider the substack $\VVV(\sigma)$ of $\XXX_{(\Sigma,\Sigma_{\nn}^0)}$,
that consists of objects
$(\pi:\SSS \to \OO_X,\alpha:\MMM\to \OO_X,\eta:\SSS\to \MMM)$
such that the condition $(**)$ holds,
where $(**)=\{$ For any geometric point $\bar{x}\to X$,
$\alpha(m)=0$ if $m\in\MMM_{\bar{x}}$ and there exists a positive
integer $n$ such that the image of $n\cdot m$ in $\overline{\MMM}_{\bar{x}}$
lies in 
$\eta(\sigma_o^{\vee}\cap\SSS_{\bar{x}})\}$.
Note that if $\Sigma$ is non-singular and $\Sigma^0_{\nn}$ is the canonical free net, i.e., $\XXX_{(\Sigma,\Sigma_{\nn}^0)}=X_{\Sigma}$, then
$(**)$ is equal to $(*)$ for any cone in $\Sigma$.
Clearly, the substack $\VVV(\sigma)$ is stable under the
torus action.
The following lemma shows that $\VVV(\sigma)$ is a closed substack
of $\XXX_{(\Sigma,\Sigma_{\nn}^0)}$.

\begin{Lemma}
\label{closed}
With the same notation as above, the condition $(**)$
is represented by a closed subscheme of $X$. In particular,
$\VVV(\sigma)\subset \XXX_{(\Sigma,\Sigma_{\nn}^0)}$ is a closed
substack.
Moreover if the base scheme $S$ is reduced, then $\VVV(\sigma)$
is reduced.
\end{Lemma}

\Proof
Let $\phi:X\to \XXX_{(\Sigma,\Sigma_{\nn}^0)}$ be
the morphism corresponding to $(\pi:\SSS \to \OO_X,\alpha:\MMM\to \OO_X,\eta:\SSS\to \MMM)$.
Since our claim is smooth local on $X$,
thus we may assume that
$\Sigma$ is a full-dimensional simplicial cone $\tau$
and $\sigma\prec \tau$.
Then by Proposition~\ref{toric2}, we have the diagram
\[
\begin{CD}
\Spec \OO_S[F]@>{p}>> \XXX_{(\tau,\tau_{\nn}^0)}@>{\pi_{(\tau,\tau_{\nn}^0)}}>>X_{\tau}=\Spec \OO_S[\tau^{\vee}\cap M]
\end{CD}
\]
where $p$ is an fppf morphism and the composite is defined by
$\tau^{\vee}\cap M\stackrel{i}{\hookrightarrow} F\stackrel{n}{\hookrightarrow} F$.
Here
$i$ is the minimal free resolution and
$n$ is defined by $e_{\rho}\mapsto n_{\rho}\cdot e_{\rho}$,
where $e_{\rho}$ is the irreducible element of $F$ corresponding
to $\rho\in \tau(1)$.
Since $p$ is fppf,
there exists \'etale locally a lifting $X\to\Spec \OO_S[F]$
of $\phi:X\to  \XXX_{(\Sigma,\Sigma_{\nn}^0)}$, and thus
we can assume $X:=\Spec \OO_S[F]$ and that 
the log structure $\MMM$ is induced by the natural map $F\to \OO_S[F]$
(cf. Proposition~\ref{AFR}).
What we have to prove is that the condition $(**)$ is represented
by a closed subscheme on $X$.
Consider the following natural diagram
\[
\xymatrix{
\tau^{\vee}\cap \sigma_{0}^{\vee}\cap M\ar[r]\ar[rd]& \tau^{\vee}\cap M\ar[d]\ar[r]^(0.6){n\circ i}& F \ar[ld]^h\\
  & \OO_S[F] &   \\
  }
\]
Set $A=\{ f\in F|\ \textup{there is}\ n\in \ZZ_{\ge 1}\ \textup{such that}\ n\cdot f\in \textup{Image}(\tau^{\vee}\cap \sigma_{0}^{\vee}\cap M)\}.$
Let $I$ be the ideal of $\OO_S[F]$ that is generated by $h(A)$.
We claim that the closed subscheme $\Spec \OO_X/I$ represents
the condition $(**)$.
To see this,
note first that
$\tau^{\vee}\cap \sigma_{0}^{\vee}\cap M$ is
$(\tau^{\vee}\cap M)\smallsetminus (\sigma^{\perp} \cap \tau^{\vee}\cap M)$
and by embeddings $(\tau^{\vee}\cap M)\subset F\subset \tau^{\vee}$,
each irreducible element of $F$ lies on a unique ray of $\tau^{\vee}$
(cf. Proposition~\ref{AFR} (2)).
Set $F=\NN^r\oplus \NN^{d-r}$ and suppose that
$\NN^{d-r}$ generates the face
$(\sigma^{\perp} \cap \tau^{\vee})$ of $\tau^{\vee}$.
Since $\tau^{\vee}\cap M$ is close to $F$, thus $I$ is the ideal generated by $\NN^r\oplus\NN^{d-r}\smallsetminus \{0\}\oplus\NN^{d-r}(\subset F)$.
Let $e_{i}$ be the $i$-th standard irreducible element of $F=\NN^d$.
Let $\bar{x}\to X=\Spec \OO_S[F]$ be a geometric point such that
$\SSS_{\bar{x}}=\gamma^{\vee}\cap M$
where $\sigma\subset \gamma\subset \tau$,
and let $\overline{\MMM}_{\bar{x}}=\langle e_{i_1},\ldots,e_{i_s}\rangle\subset F.$
In order to prove our claim, let us show that in $\OO_{X,\bar{x}}$, the
ideal $I$ coincides with the ideal
generated by
$L:=\{ m\in \overline{\MMM}_{\bar{x}}|\ \textup{there is}\ n\in \ZZ_{\ge 1}\ \textup{such that}\ n\cdot m\in \textup{Image}(\gamma^{\vee}\cap \sigma_{0}^{\vee}\cap M)\}$.
After reordering, we suppose $1 \le i_{1}<,\ldots, <i_{t}\le r<i_{t+1}<,\ldots,\le i_{s}$ ($t\le s$).
Then we have $(\OO_X/I)_{\bar{x}}=\OO_{X,\bar{x}}/(e_{i_1},\ldots,e_{i_{t}})$.
Consider the following natural commutative diagram
\[
\xymatrix@R=4mm @C=12mm{
\gamma^{\vee}\cap \sigma_0^{\vee}\cap M \ar[r] &\gamma^{\vee}\cap M \ar[d]^{=}\ar[r] &
(\gamma^{\vee}\cap M)^{\gp}\ar[r] & M\otimes_{\ZZ}\QQ \ar[r] & F^{\gp}\otimes_{\ZZ}\QQ \ar[d]\\
 & \SSS_{\bar{x}}\ar[r] & \MMM_{\bar{x}} \ar[r] & \overline{\MMM}_{\bar{x}} \ar[r] & \overline{\MMM}^{\gp}_{\bar{x}}\otimes_{\ZZ}\QQ.
}
\]
Since $\tau^{\vee}\cap \sigma_{0}^{\vee}\cap M$ is equal to $(\tau^{\vee}\cap M)\smallsetminus (\sigma^{\perp} \cap \tau^{\vee}\cap M)$ and the image of $\sigma^{\perp}(\subset M\otimes_{\ZZ}\QQ)$ in $\overline{\MMM}^{\gp}_{\bar{x}}\otimes_{\ZZ}\QQ$
is the vector space spanned by $e_{i_{t+1}},\ldots, e_{i_{s}}$,
thus
the image of $\gamma^{\vee}\cap \sigma_0^{\vee}\cap M$ in
$\overline{\MMM}_{\bar{x}}$
is contained in 
the set
\[
L':=\NN\cdot e_{i_1}\oplus\cdots \oplus \NN\cdot e_{i_{s}}\smallsetminus \NN\cdot e_{i_{t+1}}\oplus\cdots \oplus \NN\cdot e_{i_{s}}
\]
in $\overline{\MMM}_{\bar{x}}$.
Since the image of $\gamma^{\vee}\cap M$ in $\overline{\MMM}_{\bar{x}}$
is close to $\overline{\MMM}_{\bar{x}}$, we have $L=L'$.
Therefore the ideal of $\OO_{X,\bar{x}}$
generated by $e_{i_{1}},\ldots,e_{i_{t}}$
is the ideal of $\OO_{X,\bar{x}}$ generated by $L$.
Hence the
ideal $I$ coincides with the ideal
generated by $L$ in $\OO_{X,\bar{x}}$
and thus they are the same in $\OO_{X}$
because $I$ is finitely generated.
Thus $\Spec \OO_X/I$ represents
the condition $(**)$.

Finally, we will show the last assertion.
Suppose that $S$ is reduced.
With the same notation and assumption as above, $\OO_S[F]/I$ is reduced.
Since $\Spec \OO_S[F]\to  \XXX_{(\tau,\tau_{\nn}^0)}$ is fppf,
thus $\VVV(\sigma)$ is reduced.
\QED

\begin{Definition}
We shall call the closed substack $\VVV(\sigma)$
the {\it torus-invariant cycle} associated to $\sigma$.
\end{Definition}

\begin{Proposition}
Let $\tau$ and $\sigma$ be cones in a stacky fan
$(\Sigma,\Sigma_{\nn}^0)$
and suppose that $\sigma\prec\tau$.
Then we have $\VVV(\tau)\subset \VVV(\sigma)$, i.e., $\VVV(\tau)$
is a closed substack of $\VVV(\sigma)$.
\end{Proposition}

\Proof
By the definition and Proposition~\ref{closed},
it is enough to prove that for any cone $\gamma$ such that $\sigma\prec \tau\prec\gamma$, we have $\gamma^{\vee}\cap \sigma_{0}^{\vee}\cap M\subset\gamma^{\vee}\cap \tau_{0}^{\vee}\cap M$.
Clearly $\gamma^{\vee}\cap \sigma_{0}^{\vee} \subset\gamma^{\vee}\cap \tau_{0}^{\vee}$, thus our claim follows.
\QED

{\it Stabilizer groups.}
By the calculation of stabilizer group scheme
in section 3.4 and Proposition~\ref{toric1},
we can easily calculate the stabilizer group schemes
of points on a toric algebraic stack.

\begin{Proposition}
\label{ord}
Let $(\Sigma,\Sigma_{\nn}^0)$ be a stacky fan
and $\XXX_{(\Sigma,\Sigma_{\nn}^0)}$ the associated stack.
Let $\bar{x}:\Spec K\to \XXX_{(\Sigma,\Sigma_{\nn}^0)}$ be a geometric point on $\XXX_{(\Sigma,\Sigma_{\nn}^0)}$
such that $\bar{x}$ lies in $\VVV(\sigma)$,
but does not lie in any torus-invariant proper
substack $\VVV(\tau)$ \textup{(}$\sigma\prec \tau$\textup{)}. Here $K$ is an algebraically closed $S$-field.
Let $\SSS\to \OO_{K}$ be the skeleton for $\pi_{(\Sigma,\Sigma_{\nn}^0)}\circ \bar{x}:\Spec K\to \XXX_{(\Sigma,\Sigma_{\nn}^0)}\to X_{\Sigma}$.
Then $\SSS_{\bar{x}}=\sigma^{\vee}\cap M$,
and the stabilizer group scheme $G$ at $\bar{x}$
is the Cartier dual of $F^{\gp}/\iota((P)^{\gp}$.
Here $P:=(\sigma^{\vee}\cap M)/(\textup{invertible elements})$
and $\iota:P\to F$ is $t\circ r$ in Definition~\ref{toricstack} (6).
The rank of the stabilizer group scheme over $\bar{x}$, i.e.,
$\dim_{K}\Gamma(G,\OO_G)$
is
$\mult(\sigma,\sigma_{\nn}^0)=\mult (\sigma)\cdot \Pi_{\rho\in \sigma(1)}n_{\rho}$.
\end{Proposition}

\Proof
By our assumption and the definition of $\VVV(\sigma)$, clearly,
$\SSS_{\bar{x}}=\sigma^{\vee}\cap M$.
Taking account of Proposition~\ref{order} and Proposition~\ref{toric1},
the stabilizer group scheme $G$ at $\bar{x}$
is the Cartier dual of $F^{\gp}/\iota((P)^{\gp}$.
To see the last assertion, it suffices to
show that
the order of $F^{\gp}/\iota((P)^{\gp}$
is equal to $\mult (\sigma)\cdot \Pi_{\rho\in \sigma(1)}n_{\rho}$.
It follows from the fact that the multiplicity of $P$
is equal to $\mult (\sigma)$ (cf. Remark~\ref{MFRR} (3)).
\QED

\begin{Remark}
The order of stabilizer group on a point on a Deligne-Mumford
stack is an important invariant to intersection theory
with rational coefficients on it (cf. \cite{V}).

The stabilizer groups on a toric algebraic stacks
are fundamental invariants because they have data
arising from levels of the stacky fan and multiplicities
of cones in the stacky fan.
In particular, the stabilizer group scheme on a generic point on
a torus-invariant divisor $\VVV(\rho)$ ($\rho\in \Sigma(1)$)
is isomorphic to $\mu_{n_{\rho}}$.
\end{Remark}

\begin{Corollary}
\label{Deligne-Mumfordness}
Let $(\Sigma,\Sigma_{\nn}^0)$ be a stacky fan
and $\XXX_{(\Sigma,\Sigma_{\nn}^0)}$ the associated stack.
Then $\XXX_{(\Sigma,\Sigma_{\nn}^0)}$ is a Deligne-Mumford stack
if and only if $(\Sigma,\Sigma_{\nn}^0)$ is tame over $S$.
\end{Corollary}

\Proof
First of all, the ``if" direction follows from Theorem~\ref{maintoric}.
Thus we will show the
``only if" direction.
Assume that $(\Sigma,\Sigma_{\nn}^0)$ is not tame over $S$.
Then there exists a cone $\sigma\in \Sigma$
such that $r:=\mult (\sigma,\sigma_{\nn}^0)$ is not invertible
on $S$.
Then there exists an algebraically closed $S$-field $K$
such that the characteristic of $K$ is a prime divisor of $r$.
Consider a point $w:\Spec K\to \XXX_{(\Sigma,\Sigma_{\nn}^0)}\times_S K$
that corresponds to a triple $(\SSS\to \OO_{\Spec K},\MMM\to \OO_{\Spec K},\SSS\to\MMM)$
such that $\SSS=\sigma^{\vee}\cap M$.
Then by Proposition~\ref{ord},
the stabilizer group scheme of $w$
is the Cartier dual of a finite group that
has the form $\ZZ/p_1^{c_1}\ZZ\oplus \cdots \oplus\ZZ/p_l^{c_l}\ZZ$
such that $p_{i}$'s are prime numbers
and $r=\Pi_{0\le i\le l}p_i^{c_i}$.
After reordering, assume that $p_1$ is the characteristic of $K$.
Then the stabilizer group scheme of $w$ is not \'etale
over $K$. Thus by \cite[Lemma 4.2]{LM},
$\XXX_{(\Sigma,\Sigma_{\nn}^0)}$ is not a Deligne-Mumford stack.
Hence our claim follows.
\QED

\begin{Remark}
The first version \cite{Iwatoro} of this paper
only treats the case when toric algebraic stacks are Deligne-Mumford stacks,
although our method can apply to the case over arbitrary base schemes.
\end{Remark}

{\it Log structure on $\XXX_{(\Sigma,\Sigma_{\nn}^0)}$.}
A canonical log structure
on $\textup{Lis-\'et}(\XXX_{(\Sigma,\Sigma_{\nn}^0)})$
is defined as follows.
To a smooth morphism $X\to \XXX_{(\Sigma,\Sigma_{\nn}^0)}$
corresponding to a triple $(\pi:\SSS \to \OO_X,\alpha:\MMM\to \OO_X,\eta:\SSS\to \MMM)$, we attach the log structure $\alpha:\MMM\to \OO_X$.
It gives rise to a fine and saturated log structure
$\MMM_{(\Sigma,\Sigma_{\nn}^0)}$ on $\textup{Lis-\'et}(\XXX_{(\Sigma,\Sigma_{\nn}^0)})$. We shall refer it as
the {\it canonical log structure}.
Through the identification $\XXX_{(\Sigma,\Sigma_{\nn}^0)}=\XXX_{(\Sigma,T_{\Sigma},\nn)}$ (cf. Proposition~\ref{toric1}),
it is equivalent to the canonical log structure $\MMM_{(X_{\Sigma},T_{\Sigma},\nn)}$ on
$\XXX_{(\Sigma,T_{\Sigma},\nn)}$.
As in section 3.5, we have a natural morphism of log stacks
\[
(\pi_{(\Sigma,\Sigma_{\nn}^0)},\phi_{(\Sigma,\Sigma_{\nn}^0)}):(\XXX_{(\Sigma,\Sigma_{\nn}^0)},\MMM_{(\Sigma,\Sigma_{\nn}^0)})\to (X_{\Sigma},\MMM_{\Sigma})
\]
whose underlying morphism
is the canonical coarse moduli map.

\begin{Theorem}
\label{logncd}
Assume that the base scheme $S$ is regular.
Then the complement $\DDD=\XXX_{(\Sigma,\Sigma_{\nn}^0)}-T_{\Sigma}$
with reduced closed substack structure is a divisor
with normal crossings relative to $S$, and
the log structure $\MMM_{(\Sigma,\Sigma_{\nn}^0)}$ is
arising from $\DDD$.
Furthermore $\DDD$ coincides with the union of
$\cup_{\rho\in\Sigma(1)}\VVV(\rho)$.
\end{Theorem}

For the proof we need the following lemma.

\begin{Lemma}[log EGA IV 17.7.7]
\label{log17.7.7}
Let $(f,\phi):(X,\LLL)\to (Y,\MMM)$
and $(g,\psi):(Y,\MMM)\to (Z,\NNN)$ be
morphisms of fine log schemes respectively.
Suppose that $(g,\psi)\circ (f,\phi)$ is log smooth
and $g$ is locally of finite presentation.
Assume further that $(f,\phi)$ is a strict faithfully flat
morphism that is locally of finite presentation.
Then $(g,\psi)$ is log smooth.
\end{Lemma}

\Proof
Let
\[
\xymatrix{
(T_0,\PPP_0)\ar[r]^s \ar[d]^i &  (Y,\MMM)\ar[d] \\
(T,\PPP) \ar[r]^t & (Z,\NNN) \\
}
\]
be a commutative diagram of fine log schemes
where $i$ is a strict closed immersion defined by a nilpotent
ideal $I\subset \OO_T$.
It suffices to show that there exists \'etale locally on $T$
a morphism $(T,\PPP)\to (Y,\MMM)$ filling in the diagram.
To see this, after replacing $T_0$ by an \'etale cover
we may suppose that $s$ has a lifting
$s':(T_0,\PPP_0)\to (X,\LLL)$.
Since $(g,\psi)\circ (f,\phi)$ is log smooth,
there exists \'etale locally on $T$ a lifting $(T,\PPP)\to (X,\LLL)$ of $s'$
filling in the diagram. Hence our assertion easily follows.
\QED

{\it Proof of Theorem~\ref{logncd}}.
By Proposition~\ref{toric1} and Proposition~\ref{localalg} (together with
Remark~\ref{fullgene}),
there is an fppf strict morphism
$(X,\LLL)\to (\XXX_{(\Sigma,\Sigma_{\nn}^0)},\MMM_{(\Sigma,\Sigma_{\nn}^0)})$
from a scheme $X$ such that
the composite $(X,\LLL)\to (S,\OO_S^*)$ is log smooth.
Here $\OO_S^*$ denotes the trivial log structure on $S$.
Then by applying Lemma~\ref{log17.7.7}
we see that $(\XXX_{(\Sigma,\Sigma_{\nn}^0)},\MMM_{(\Sigma,\Sigma_{\nn}^0)})
\to (S,\OO_S^*)$ is log smooth.
Let $U\to \XXX_{(\Sigma,\Sigma_{\nn}^0)}$ be a smooth surjective morphism
from an $S$-scheme $U$.
Then by \cite[Theorem 4.8]{FK}
there exists \'etale locally on $U$
a smooth morphism $U\to S\times_{\ZZ}\ZZ[P]$
where $P$ is a toric monoid
(\cite[Theorem 4.8]{FK} worked over a field, but it
is also applicable to our case).
In addition, the natural map $P\hookrightarrow \OO_S[P]$ gives rise to
a chart for $\MMM_{(\Sigma,\Sigma_{\nn}^0)}$.
Since $\XXX_{(\Sigma,\Sigma_{\nn}^0)}$ is smooth over $S$, thus
after shrinking $S\times_{\ZZ}\ZZ[P]$
we may suppose that $P$ is a free monoid.
Note that the support of $\overline{\MMM}_{(\Sigma,\Sigma_{\nn}^0)}$
is the complement $\XXX_{(\Sigma,\Sigma_{\nn}^0)}-T_{\Sigma}$.
Therefore $\XXX_{(\Sigma,\Sigma_{\nn}^0)}-T_{\Sigma}$
with reduced substack structure
is a divisor with normal crossings relative to $S$.
Next we will prove that
$\cup_{\rho\in\Sigma(1)}\VVV(\rho)$ is reduced. (Set-theoretically $\cup_{\rho\in\Sigma(1)}\VVV(\rho)=\XXX_{(\Sigma,\Sigma_{\nn}^0)}-T_{\Sigma}$.)
It follows from Lemma~\ref{closed}.
Finally, we will show $\MMM_{(\Sigma,\Sigma_{\nn}^0)}\cong \MMM_{\DDD}$
where $\MMM_{\DDD}$ is the log structure associated to $\DDD$.
Note that $\MMM_{\DDD}=\OO_{\XXX_{(\Sigma,\Sigma_{\nn}^0)}}\cap i_{(\Sigma,\Sigma_{\nn}^0)*}\OO_{T_{\Sigma}}^*$ where $i_{(\Sigma,\Sigma_{\nn}^0)}:T_{\Sigma}\to \XXX_{(\Sigma,\Sigma_{\nn}^0)}$ is the canonical torus embedding.
Hence $\MMM_{\DDD}$ is a {\it subsheaf} of $\OO_{\XXX_{(\Sigma,\Sigma_{\nn}^0)}}$. Thus we may and will assume that $\XXX_{(\Sigma,\Sigma_{\nn}^0)}=\XXX_{(\tau,\tau_{\nn}^0)}$ where $\tau$ is a full-dimensional cone.
Let $p:\Spec \OO_S[F]\to \XXX_{(\tau,\tau_{\nn}^0)}$ be an fppf cover
given in Proposition~\ref{toric2}.
By the fppf descent theory for fine log structures (\cite[Corollary A.5]{OL})
together with the fact that $\MMM_{\DDD}$ is a subsheaf $\OO_{\XXX_{(\tau,\tau_{\nn}^0)}}$, it is enough to show $p^*\MMM_{(\tau,\tau_{\nn}^0)}\cong p^*\MMM_{\DDD}$. Since $p^*\MMM_{(\tau,\tau_{\nn}^0)}$
is induced by the natural map $F\to \OO_S[F]$ (cf. Proposition~\ref{toric1} and
Proposition~\ref{ex2}), thus we have $p^*\MMM_{(\tau,\tau_{\nn}^0)}\cong p^*\MMM_{\DDD}$.
\QED

\begin{Proposition}
\label{logetale}
With the same notation and assumptions as above,
suppose further that $(\Sigma,\Sigma_{\nn}^0)$ is tame over $S$.
Then $(\pi_{(\Sigma,\Sigma_{\nn}^0)},\phi_{(\Sigma,\Sigma_{\nn}^0)}):(\XXX_{(\Sigma,\Sigma_{\nn}^0)},\MMM_{(\Sigma,\Sigma_{\nn}^0)})\to (X_{\Sigma},\MMM_{\Sigma})$ is {\it Kummer log \'etale}.
In particular, there exists an isomorphism of $\OO_{\XXX_{(\Sigma,\Sigma_{\nn}^0)}}$-modules
\[
\Omega^1((\XXX_{(\Sigma,\Sigma_{\nn}^0)},\MMM_{(\Sigma,\Sigma_{\nn}^0)}))/(S,\OO_S^*))\stackrel{\sim}{\to}\OO_{\XXX_{(\Sigma,\Sigma_{\nn}^0)}}\otimes_{\ZZ}M,
\]
where $\Omega^1((\XXX_{(\Sigma,\Sigma_{\nn}^0)},\MMM_{(\Sigma,\Sigma_{\nn}^0)}))/(S,\OO_S^*))$ is the sheaf of log differentials (cf. \cite[5.6]{FK}).
\end{Proposition}

\Proof
It follows from the next Lemma, Proposition~\ref{toric1} and Proposition~\ref{localalg}.
\QED

\begin{Lemma}
\label{kummerlog}
With the same notation as in Proposition~\ref{ex2},
if the order of $F^{\gp}/\iota^{\gp}(P^{\gp})$ is invertible
on $R$, then
\[
(f,\phi):(\Spec R[F],\MMM_F)\to (\Spec R[P],\MMM_P)
\]
is Kummer log \'etale.
In particular, the induced morphism (cf. section 3.2)
\[
([\Spec R[F]/G],\MMM)\to (\Spec R[P],\MMM_P)
\]
is Kummer log \'etale. Here $G$
is the Cartier dual of $F^{\gp}/\iota(P^{\gp})$
and $\MMM$
is the log structure associated to a natural chart $F\to R[F]$.
\end{Lemma}

\Proof
Since $(f,\phi)$ is an admissible FR morphism, thus $(f,\phi)$
is Kummer.
By the toroidal characterization of log \'etaleness
(\cite[Theorem 3.5]{log}),
$(f,\phi):(\Spec R[F],\MMM_F)\to (\Spec R[P],\MMM_P)$ is log \'etale.
Since $\Spec R[F]\to [\Spec R[F]/G]$ is \'etale, thus our claim follows.
\QED

{\it Proof of Proposition~\ref{log1}.}
With the same notation as in Proposition~\ref{log1}, we may assume that $X$ is an affine simplicial toric variety.
Then Proposition~\ref{log1}
follows from Proposition~\ref{toric1}, Theorem~\ref{logncd} and Proposition~\ref{logetale}.
\QED

\end{document}